\documentclass{amsart}[12pt]
\usepackage{amsmath,amssymb}
\usepackage[all]{xy}

\newcommand{\C}{{\mathbb{C}}}
\newcommand{\N}{{\mathbb{N}}}

\newcommand{\R}{{\mathbb{R}}}

\newcommand{\Ch}{{\mathcal C}}

\newcommand{\Fh}{{\mathcal F}}

\newcommand{\Uh}{{\mathcal U}}

\newcommand{\Zh}{{\mathcal Z}}

\newcommand{\F}{\mathcal{F}}
\newcommand{\G}{\mathcal{G}}
\newcommand{\A}{A}
\newcommand{\K}{\mathcal{K}}
\newcommand{\Mul}{\mathcal{M}}
\renewcommand{\P}{\mathcal{P}} 
\newcommand{\M}{ \mathbb{M}}

\newcommand{\AF}{\mathrm{AF}}
\newcommand{\Aff}{\mathrm{Aff}}
\newcommand{\AH}{\mathrm{AH}}
\newcommand{\ASH}{\mathrm{ASH}}
\newcommand{\AI}{\mathrm{AI}}

\newcommand{\be}{\mathbf{1}}

\newcommand{\card}{\mathrm{card}}
\newcommand{\diam}{\mathrm{diam}}
\newcommand{\dist}{\mathrm{dist}}

\newcommand{\ev}{\mathrm{ev }}

\newcommand{\halb}{\frac{1}{2}}
\newcommand{\her}{\mathrm{her}}

\newcommand{\Hom}{\mathrm{Hom}}
\newcommand{\id}{\mathrm{id}}
\newcommand{\kker}{\mathrm{ker}}
\newcommand{\Kk}{\mathrm{K}}
\newcommand{\KK}{\mathrm{KK}}
\newcommand{\KL}{\mathrm{KL}}
\newcommand{\mmin}{\mathrm{min}}

\newcommand{\TT}{\mathrm{T}}
\newcommand{\TAF}{\mathrm{TAF}}
\newcommand{\tr}{\mathrm{tr}}
\newcommand{\UHF}{\mathrm{UHF}}

\newcounter{number}[section]

\newenvironment{nummer}{\refstepcounter{number}{\noindent\arabic{section}.\arabic{number}}}{}

\newcommand{\bn}{\noindent \begin{nummer} \rm}
\newcommand{\en}{\end{nummer}}

\newenvironment{thm}{ \noindent {\sc Theorem:} \it}{}
\newenvironment{lms}{\noindent {\sc Lemma:} \it}{}
\newenvironment{lem}{\noindent {\sc Lemma:} \it}{}
\newenvironment{props}{\noindent {\sc Proposition:} \it}{}
\newenvironment{prop}{\noindent {\sc Proposition:} \it}{}

\newenvironment{df}{\noindent {\sc Definition:} \it}{}

\newenvironment{nots}{\noindent {\sc Notation:} }{}

\newenvironment{nproof}{\noindent {\sc Proof:}}{\mbox{}\hfill 
\rule[-.2ex]{.25em}{1.8ex}}

\begin{document}

\title[Commutative subalgebras of simple 
$C^{*}$-algebras]{Commutative 
$C^{*}$-subalgebras of simple stably finite
$C^{*}$-algebras with real rank zero}

\author{Ping Wong Ng}
\address{Department of Mathematics\\
University of Louisiana\\
Lafayette, LA\\
70504\\ 
USA}
\email{png@louisiana.edu}

\author{Wilhelm Winter}
\address{Mathematisches Institut der Universit\"at M\"unster\\
Einsteinstr. 62\\ D-48149 M\"unster}

\email{wwinter@math.uni-muenster.de}

\date{\today}
\thanks{{\it Supported by:} DFG (through the SFB 478), EU-Network  Quantum Spaces - Noncommutative 
\indent Geometry (Contract No. HPRN-CT-2002-00280)}

\setcounter{section}{-1}

\begin{abstract}
Let $X$ be a second countable, path connected, compact metric space
and let $A$ be a unital separable simple
nuclear $\mathcal{Z}$-stable real rank zero $C^*$-algebra.
We classify all the unital $*$-embeddings (up to approximate unitary 
equivalence)
of $C(X)$ into $A$.
Specifically, we provide an existence and a uniqueness theorem for
unital $*$-embeddings from $C(X)$ into $A$.   
\end{abstract}

\maketitle

\section{Introduction}

    In this paper, we study the classification problem for
embeddings of a  
commutative $C^*$-algebra $C(X)$ into a simple, real rank zero,
stably finite $C^*$-algebra 
$A$ (with other properties).  
This problem is closely related to, though independent of, the 
classification problem for the simple $C^*$-algebra $A$ itself.
Moreover, techniques from either subject often carry over to the other and
have interesting implications (for other areas as well!).  

    A complete solution to 
the classification problem would consist of two parts: an existence theorem
and a uniqueness theorem. 
One starting point for the existence theorem is  
a result of Pimsner \cite{Pimsner} which shows that 
there is an interesting $*$-monomorphism $\phi : C(T^2) \rightarrow A$ of
$C(T^2)$ into a unique simple $\AF$-algebra $A$ with 
$\Kk_0(A) \cong \mathbb{Z}[1/2] \oplus \mathbb{Z}$ with lexicographic order
and order unit $(1,0)$ (here $T^2$ is the $2$-torus).  
This $*$-monomorphism actually
induces an isomorphism on the rational $\Kk_0$ groups.  
The example of Pimsner exhibits the higher dimensional features of the 
$\AF$-algebra which is also, at the same time, a noncommutative zero 
dimensional space (of course, every commutative $C^*$-algebra can
be embedded in a commutative $C^*$-algebra with spectrum being a 
Cantor set - but not so with interesting induced map in $\mathrm{K}$-theory). 
This phenomenon was subsequently 
intensively studied by Dadaralat, Elliott and 
Loring, who showed that many group homomorphisms from $\Kk_0(C(X))$
to $\Kk_0(A)$ can be realized by $C^*$-homomorphisms from $C(X)$
to $A$.  More specifically, 
in \cite{EL}, Elliott and Loring showed that for a unital simple
$\AF$-algebra $A$, a group homomorphism 
$\mu : \Kk_0(C(T^2)) \rightarrow \Kk_0(A)$ can be realized by a 
unital $*$-homomorphism
$\phi : C(T^2) \rightarrow A$ if and only if $\mu$ preserves
the
order unit and sends every element of the reduced $\Kk_0$ group of $C(T^2)$
to an element of $\kker (\tau_*)$ for every unital trace $\tau \in \TT (\A)$.
In \cite{DadL},  Dadarlat and Loring showed that given a finite
CW complex $X$, there are a unital $\AF$-algebra $A$ and a unital
$*$-monomorphism $\phi : C(X) \rightarrow A$ such that $\phi$
induces an isomorphism on the rational $\Kk_0$ groups.
In particular, if $X$ is connected then $A$ can be taken to be the unique
unital simple $\AF$-algebra with $\Kk_0$ group $\mathbb{Z}[1/2] \oplus 
\mathbb{Z}^n$, where $n$ is the rank 
of the reduced $\Kk_0$ group of $C(X)$.  Here,
the positive cone is $\{ (r, m) : r \in \mathbb{Z}[1/2], \makebox{ }
r > 0 \makebox{  and  } m \in 
\mathbb{Z}^n \}$ and the order unit is $(1,0)$. 
In their paper \cite{DadL}, 
Dadarlat and Loring also proved other existence 
results.  For instance, let $X$ be a finite CW-complex and let $A$
be the tensor product of a unital simple $\AF$-algebra with the $\UHF$-algebra
with dimension group being the rational numbers.   
If $\mu : \Kk_0( C(X)) \rightarrow \Kk_0 (A)$ is a group homomorphism such that
(i)  $\mu ([\be])  = [\be]$ and (ii) $\mu$ brings the reduced $\Kk_0$ group
of $C(X)$ to $\bigcap_{\tau \in \TT (A) } \kker ( \tau_*)$, then there
is a unital $*$-homomorphism $\phi : C(X) \rightarrow A$ such that 
$\Kk_0 (\phi) = \mu$ (this generalizes the result of Elliott and Loring).

     Li and Lin have generalized the results of Dadarlat, Elliott, Loring
and Pimsner to a large class of codomain algebras $A$.
In \cite{Liexistence},  Li showed that
given a finite CW complex $X$, given a unital simple $\AH$-algebra $A$ 
with real
rank zero and bounded dimension growth and given $\alpha \in \KK(C(X), A)$,
$\alpha$ can be realized by a unital $*$-homomorphism $\phi : C(X) \rightarrow
A$ if and only if $\alpha \in \KK( C(X), A)_{+, D(A)}$ ($\KK(C(X), A)_{+, D(A)}$
is a subset of the $\KK$-group $\KK(C(X), A)$;  see \cite{Liexistence} for
the precise definition).  
In \cite{Linembedding}, Lin showed that the codomain algebra $A$ (in Li's
result) can be replaced by an arbitrary unital simple separable nuclear  
$C^*$-algebra with real rank zero, stable rank one and weak unperforation.
We note that the results of Dadarlat, Elliott, Li, Lin, Loring 
and Pimsner were proven after (and in the case of $\AF$-algebras, long after!)
the classification result for the corresponding simple codomain algebras.
We also note that
their results are incomplete in the sense that they do not fully take
into account the tracial simplex (an important classification invariant).
The best result in this respect 
was that of Li (in \cite{Liexistence}) where for 
$A$ having unique trace (in addition to other already mentioned conditions),
the $*$-homomorphism can be constructed to realize any strictly positive
probability measure on $X$.   
In this paper, we complete the existence theorem by fully taking into 
account the tracial simplex.  Our result is as follows:

\bn
\begin{thm}
Let $X$ be a second countable, path connected, compact metric space
and let $A$ be a unital separable simple
nuclear $C^*$-algebra with real rank zero, stable rank one and weak 
unperforation.
Let $\alpha \in \KL (C(X), A)_{+, 1}$ and let 
$\lambda : \TT (A) \rightarrow \TT (C(X))$ be an affine continuous map such 
that 
\begin{enumerate}
\item[(a)] if $h \in \Aff(\TT (C(X)))$ with $h \geq 0$ and $h$ is not the zero
function then $\Aff(\lambda)(h)(\tau) > 0$ for all $\tau \in \TT (A)$; and
\item[(b)]  for every projection $p \in C(X) \otimes \K$,
$\lambda(\tau)(p) = \tau_{\ast}(\alpha([p]))$ for all $\tau \in \TT (A)$.
\end{enumerate}
Then there is a unital $*$-monomorphism $\phi : C(X) \rightarrow A$ such 
that $\KL(\phi) = \alpha$ and $\TT (\phi) = \lambda$.
\label{thm:MainExistence}
\end{thm}
\en

     Another part of the classification problem (for embeddings of $C(X)$
into simple real rank zero stably finite etc. $A$) is the uniqueness theorem.
One early result, without assuming that the codomain algebra is stably 
finite, is the result in \cite{BDF} where Brown, Douglas and Fillmore
showed that for two unital $*$-monomorphisms $\phi, \psi : C(X) \rightarrow
\mathbb{B}(H)/ \K(H)$,  $\phi$ and $\psi$ are unitarily equivalent if and only
if $\KK(\phi) = \KK(\psi)$. 
Another early result is in \cite{GongLin} where Gong and Lin showed that 
for $X$ a compact metric space and $A$ a unital simple separable nuclear 
$C^*$-algebra with real rank zero, stable rank one, weak unperforation and
unique tracial state,
two unital 
$*$-monomorphisms $\phi, \psi : C(X) \rightarrow A$ are approximately 
unitarily equivalent if and only if $\KL (\phi) = \KL (\psi)$ and 
$\tau \circ \phi = \tau \circ \psi$ for all $\tau \in \TT (A)$.
It is not hard to replace the unique tracial state condition (in Gong and
Lin's result) by the 
(slightly) weaker condition of countably many extreme tracial states.
However, moving beyond these conditions of having small tracial simplex
is difficult.  This is closely related to certain recent classification
problems, for simple $C^*$-algebras, 
where a major difficulty lies in replacing the small tracial
simplex condition by arbitrary tracial simplexes (see, for example 
\cite{NBrown}, \cite{Lintraces}, \cite{Ng}, \cite{Wi5} and \cite{Wi6}). 

      Recently in \cite{Linhomomorphism}, 
Lin replaced the small tracial simplex condition in the 
codomain algebra $A$ by the additional assumption that $A$ had 
tracial rank zero (i.e., the $\TAF$ property;  this is a strong property
which, together with the assumptions of nuclearity and universal coefficient
theorem, imply that the algebra involved is a simple unital $\AH$-algebra
with bounded dimension growth and real rank zero).
We note that Lin used this result (in \cite{Linhomomorphism}) to 
generalize an interesting theorem of Kishimoto.  Specifically, Lin showed 
that if $A$ is a simple unital $\AH$-algebra with real rank zero and 
bounded dimension growth, and if $\alpha$ is an approximately inner 
$*$-automorphism of $A$ with the tracial Rokhlin property, then
the crossed product $A \times_{\alpha} \mathbb{Z}$ is a unital simple 
$\AH$-algebra with bounded dimension growth and real rank zero.

      In this paper, we show that the restriction on the tracial simplex
of $A$ can be removed without assuming that $A$ is $\TAF$.  However, we
need to assume that $A$ is $\mathcal{Z}$-stable where $\mathcal{Z}$ is 
the Jiang-Su algebra.  This 
condition is implied by approximate divisibility
and is currently of interest in the
classification program for simple $C^*$-algebras (more below).
Our uniqueness result is as follows:

\bn
\begin{thm}
Let $X$ be a second countable, path connected, compact metric space
and let $A$ be a unital separable simple
nuclear real rank zero $\mathcal{Z}$-stable $C^*$-algebra.
Let $\phi, \psi : C(X) \rightarrow A$ be unital $*$-monomorphisms.
Then $\phi$ and $\psi$ are approximately unitarily equivalent if and only if
$\KL (\phi) = \KL (\psi)$ and $\tau \circ \phi = \tau \circ \psi$ for all
$\tau \in \TT (A)$.
\label{thm:MainUniqueness}
\end{thm}
\en
  
    The \emph{Jiang-Su algebra} $\Zh$ is the unique unital simple nontype I
$\ASH$-algebra (i.e., approximately subhomogeneous $C^*$-algebra)  
such that (a) $\Zh$ has
the same invariant as the complex numbers and (b) $\Zh$ is strongly
self-absorbing (see \cite{TW1}, \cite{TW2}, \cite{TW3}, \cite{DadT}, and
\cite{Wi7}; see also \cite{RW1} for different descriptions of $\Zh$).  
A $C^*$-algebra $A$ is said to be \emph{$\Zh$-stable} if $A \otimes \Zh 
\cong A$.  $\Zh$-stability implies many nice properties
from the point of view of classification theory.  Among other things,
if $A$ is a unital separable simple $\Zh$-stable $C^*$-algebra then $A$ is
either purely infinite or stably finite; and if, in addition,  $A$ is 
stably finite then $A$ must have stable rank one and weak unperforation
(see \cite{GJS}, \cite{JS1} and \cite{R4}).
Many simple $C^*$-algebras are $\Zh$-stable - including all the ones that
have been classified using $\mathrm{K}$-theory invariants (see \cite{TW2});  
and $\Zh$-stability
has played an interesting role in recent advances in classification theory.
(See, for example, \cite{Wi5}, \cite{Wi6}, \cite{Wi7} and \cite{LinApp}.)  
Finally, it has been put forward in several places that 
$\Zh$ should play the role in general classification theory that $\mathcal{O}_{\infty}$
plays for the classification of simple nuclear purely infinite $C^*$-algebras. 

    We note that the subject matter of this paper is intimately
related to the question of when a unital $*$-homomorphism 
$\phi : C(X) \rightarrow A$ can be pointwise approximated by finite dimensional
$*$-homomorphisms. (See, for example, \cite{GongLin} and 
\cite{Finitedimensional} 
for more details and references.)  
This
is closely related to Lin's interesting result that pairs of  
almost commuting self-adjoint matrices are uniformly close to pairs of
exactly commuting self-adjoint matrices (see, for example, \cite{Lincommute}).
We also note that an interesting problem would be to replace the commutative
algebra $C(X)$ (in our results) by an 
arbitrary nuclear residually finite 
dimensional $C^*$-algebra.  For instance, results  of this type would
lead to a proof that every simple unital separable nuclear quasidiagonal
$C^*$-algebra is $\AF$-embeddable (this is not even known for inductive
limits of residually finite dimensional type I $C^*$-algebras).  
This in itself is interesting, but
it may also in turn lead to a proof
of appropriate tracial approximation properties and perhaps classification
(see, for example, \cite{NBrown} and \cite{Lintraces}).

   Also, Theorem \ref{thm:MainUniqueness} is not known (except for
special cases) when the codomain algebra is not real rank zero.
Indeed, we do not even have a full uniqueness theorem for the
case where the codomain algebra is a real rank one simple unital 
$\AH$-algebra with bounded dimension growth.  Hence, for this class
of algebras, even though there is a complete classification using
$\mathrm{K}$-theory invariants, the general theory is still incomplete.  
One consequence of this defect is that the structure 
of the automorphism
groups of these algebras is still not known.  (See \cite{EGL},
\cite{LinTR1}, \cite{NgRuiz1}, \cite{NgRuiz2}, \cite{NT},
\cite{Mygind1}, \cite{Mygind2}).      

   Finally, the results of this paper are connected to but distinct
from the existence and uniqueness theorems found in 
\cite{Dad}, \cite{DE} and \cite{LinUniqueness}.  The theorems
in these papers require that the domain algebra be a simple $\TAF$
$C^*$-algebra.   
The first immediate result is that the tracial simplex need not 
be taken into account (whereas for commutative domain algebras, 
the tracial simplex is needed and results in much extra work).  
Secondly, that the domain algebra is simple $\TAF$ leads to simplication
and special techniques for uniqueness, which do not work for 
commutative domain algebras - this is why, for example, Lin had 
to prove a separate uniqueness theorem for commutative domain algebras
in \cite{Linhomomorphism}.  Finally, since all the above 
existence theorems require uniqueness, similar remarks hold for the
existence theorems.  

   In this paper, we will use ideas from $\mathrm{K}$-theory and $\KK$-theory.
A good reference for this is \cite{Linbook}.
We will also use the notions of \emph{finite decomposition rank} and
$\Zh$-\emph{stability} which have recently played important roles
in classification theory.  Good references for this are
\cite{KW} and \cite{TW2}.

    Throughout this paper, ``c.p.c." will denote
``completely positive contractive".  
For $C^*$-algebras $A, B$, $\underline{\Kk}(A)$ is the total $\mathrm{K}$-theory of $A$, 
$\mathbb{P}(A)$ is the class of projections defined in \cite{Linbook}
6.1.1 (the equivalence classes of the projections in 
$\mathbb{P}(A)$ generate $\underline{\Kk}(A)$ as a group),
$\KK(A, B)$ is the $\KK$-group of $A$ and $B$, and (when either $A$ or $B$
satisfies the UCT) 
$\KL (A, B) = \Hom_{\Lambda}(\underline{\Kk}(A), \underline{\Kk}(B))$ is the 
$\KL $ group of $A$ and $B$ (here, $\Lambda$ is the collection of 
Bockstein operations).  
We also note that if $\phi : A \rightarrow B$ is a contractive 
completely positive
map which is ``sufficiently almost multiplicative" then 
$\phi$ gives a well-defined map from a finite subset of 
$\mathbb{P}(A)$ into $\underline{\Kk}(B)$;  
a discussion of this can be found in \cite{LinTR0} page 8, and we will
be using the notation contained there.  
For more details on all of the above, see \cite{Linbook} and 
\cite{LinTR0}. 

\section{The existence theorem}

\bn
In this section, we prove an existence theorem for $*$-monomorphisms
from $C(X)$ into a simple real rank zero $C^*$-algebra $\A$ (with 
appropriate additional properties).
We first need a lemma concerning 
c.p.c. almost multiplicative maps from $C(X)$ into a matrix algebra
(cf. \cite{Linbook} Lemma 6.2.7):

\begin{lem}   
     Let $X$ be a compact metric space.  Let $\epsilon 
 > 0$ and a finite subset $\F \subseteq C(X)$ be given.\\ 
\noindent Then there 
exists $\delta > 0$ and a finite subset $\G \subseteq C(X)$ satisfying
the following:   

  For every positive integer $n \geq 1$, if $L : C(X) \rightarrow 
\M_n( \mathbb{C})$
is a $\G - \delta$ multiplicative c.p.c.  map   
then there exists a $*$-homomorphism $h : C(X) \rightarrow \M_n$ with
$h(\be_{C(X)}) = r$ such that 
\[
 \tr (\be - r) < \epsilon
\]
and 
\[
\| L(f) - (\be - p) L(f) (\be - p) + h(f) \| < \epsilon,
\]
for all $f \in \F$.\\ 
(In the above, $\tr$ is the unique normalized trace on $\M_n$.)
\label{lem:LinbookMatrix}
\end{lem}
\en

\bn
   Next, we will also need the following theorem concerning the approximation
of unital positive linear maps between tracial
state spaces (of commutative $C^*$-algebras) by
convex combinations of maps coming from $*$-homomorphisms, 
which is due to Li \cite{Litraces}:

\begin{thm}  Let $X$ be a path connected compact metric space. Let a  
finite subset $\F \subseteq \Aff (\TT (C(X)))$ and $\epsilon > 0$ be given.\\
\noindent Then there is 
an $N > 0$ with the following property:

  For any unital positive linear map $\xi: \Aff (\TT (C(X)))  \rightarrow 
\Aff (\TT (C(Y)))$, where $Y$ is an arbitrary compact metric space, there 
are $N$ $*$-homomorphisms 

\[
\phi_1, \phi_2, ..., \phi_N : C(X) \rightarrow C(Y)
\]
such that each $\phi_i$ is homotopy equivalent to a point evaluation (i.e.,
$\phi_i$ is homotopy equivalent to a $*$-homomorphism of the form 
$C(X) \rightarrow C(Y) : f \mapsto f(x_0) \be_{C(Y)}$ for some point
$x_0 \in X$) and such that   
\[
| \xi(f)(\tau) - (1/N) \sum_{i=1}^N \tau(\phi_i (f)) | < \epsilon, 
\forall f \in \F, \forall \tau \in T(C(Y)).
\]
\label{thm:Litraces}
\end{thm}
\en

\bn
To continue, we fix some notation:  Let $X$ be a compact second countable 
metric space. Let $A$ be a unital $C^*$-algebra.
We denote by $\KL (C(X), A)_{+, 1}$ the set of all $\alpha \in  
\KL (C(X), A)$ which satisfy the following two conditions:
\begin{enumerate}
\item[i.] $\alpha([\be_{C(X)}]) = [ \be_{A}] $ (for the induced map between
$\Kk_0$ groups).
\item[ii.] 
$\alpha( \Kk_0 (C(X))_+ - \{ 0 \} ) \subseteq \Kk_0 (A)_+ - \{ 0 \}$.   
\end{enumerate}
\en

\bn
   Next, we need the following result which follows from the $\TAF$-property
and \cite{Linbook} Theorem 6.1.11. (The proof can be found in 
\cite{Linbook} 
Lemma 6.2.8.  Note that the proof of this Lemma does not really
require that the spectrum of the domain algebra be a finite CW complex.
See also \cite{Dad} Theorem 5.6 and Corollary 5.7.) 

\begin{thm}  Let $A$ be a unital simple nontype I $\AH$-algebra with bounded
dimension growth and real rank zero.
Suppose that $X$ is a compact, second countable, path connected metric space,
$\P \subseteq \mathbb{P}(C(X))$ is a finite subset, and  
\[
\alpha \in \KL ( C(X), A)_{+,1}.
\]
\noindent Then there exists a sequence of unital c.p.c. linear maps 
\[
\phi_n : C(X) \rightarrow A
\]
\noindent such 
that
\[
[ \phi_n ]|\P  = \alpha | \P
\]
and
\[
\| \phi_n(fg) - \phi_n(f) \phi_n(g) \| \rightarrow 0
\]
as $n \rightarrow \infty$, for all $f, g \in C(X)$.   
\label{thm:Linexistence}
\end{thm} 

(For the definition of $[\phi_n] | \P$ we refer the reader to 
\cite{LinTR0} page 8.  As mentioned in the last paragraph of the
introduction, we will be
using the notation contained there.)

\begin{nproof}[Sketch of proof]

     As already indicated, the proof can be found in 
\cite{Linbook}  Lemma 6.2.8.  Here, we sketch the proof for the 
case where $C(X)$ has torsion-free $K_0$ and $K_1$ groups.

     By \cite{Linbook} Theorem 6.1.11, let $\phi_1 : C(X) \rightarrow
\M_n(A)$ be an almost multiplicative map  and 
let $\phi_2 : C(X) \rightarrow \M_{n - 1}(A)$ be a unital $*$-homomorphism
with finite dimensional range such that 
$$[\phi_1] | \P = \alpha | \P + [ \phi_2] | \P$$
(Of course, we need $\phi_1$ to be almost multiplicative on a sufficiently
large finite subset of $C(X)$ and for a sufficiently small positive
real number.)

     Since $\M_n(A)$ is tracially $\AF$,   
we can find a projection $p \in \M_{n}(A)$         
such that $p$ is Murray-von Neumann equivalent to 
a subprojection of $1_{A}$, and we can find
a unital c.p.c. almost multiplicative map
$\phi_3 : C(X) \rightarrow p \M_n (A) p$ and 
a unital $*$-homomorphism $\phi_4 : 
C(X) \rightarrow (1 -p) \M_n ( A) (1 - p)$ 
such that $\phi_1$ can be approximated by $\phi_3 \oplus \phi_4$
in norm on a sufficiently large finite subset of $C(X)$, and such that
$$[\phi_1] | \P = [\phi_3] | \P + [ \phi_4] | \P$$ 
(Of course, we are assuming that all the maps involved are sufficiently
multiplicative so that all the above expressions are well-defined etc.)

    By conjugating with a unitary if necessary, we may assume that 
$p \leq 1_{A}$.
Let $x_0 \in X$ be an arbitrary point. 
Then take $\phi : C(X) \rightarrow \Mul$ to be the unital c.p.c. map 
given by $\phi(f) =_{df} \phi_3(f) + f(x_0) (1_{A} - p)$, for all
$f \in C(X)$.      
Assuming that all the maps involved were chosen to be sufficiently 
multiplicative, we can check that 
$$[\phi] | \P = \alpha | \P$$
(This last equality is nontrivial!)
\end{nproof}

\en
    
Finally, towards the existence theorem, we need the following approximate
uniqueness result (for tracial rank zero codomains) of Lin (which can
be found in \cite{Linhomomorphism} Theorem 4.6;  see also 
\cite{Finitedimensional}):

\bn
\begin{thm} Let $X$ be a compact metric space, $\epsilon > 0$ and $\F 
\subseteq C(X)$ be a finite subset.  
Let $\nu > 0$ be such that $| f(x) - f(y) | < \epsilon/8$ if $d(x,y) < \nu$
for all $f \in \F$ and all $x, y \in X$ ($d$ is the metric on $X$).   
Then, for any integer $s \geq 1$, any
finite $\nu/2$-dense subset $\{ x_1, x_2,..., x_m \}$ of $X$ for which
$O_i \cap O_j = \emptyset$ for $i \neq j$, where  
$O_i =_{df} \{ x \in X : \dist (x, x_i) < \nu/(2s) \}$ 
and any $1/(2s) > \sigma > 0$, there exists $\gamma > 0$, a finite    
subset $\G \subseteq C(X)$, $\delta > 0$ and a finite subset $\P \subseteq
\mathbb{P}(C(X))$ satisfying the following:

    For any unital separable simple nuclear
$C^*$-algebra $A$ with tracial rank zero,
any $\G - \delta$-multiplicative unital c.p.c. 
linear maps  $\phi, \psi : C(X) \rightarrow A$ with 
$\tau \circ \phi(g)$ within $\gamma$ of $\tau \circ \psi(g)$ 
for all $g \in \G$,  if 
\begin{enumerate}
\item[(i)] $\mu_{ \tau \circ \phi} (O_i), \mu_{\tau \circ \psi} (O_i) > 
\sigma \nu$ for all $i$ and for all $\tau \in \TT (A)$, and
\item[(ii)] $[ \phi ] | \P = [ \psi ]| \P$
\end{enumerate}
then there exists a unitary $u \in A$ such that 
$u \phi (f) u^*$ is within $\epsilon$ of $\psi(f)$ for all $f \in \F$.  

    Finally, if in the above, the elements of $\F$ all have norm less than
or equal to one, then we can choose $\G$ so that its elements all have norm
less than or equal to one. 
\label{thm:Linapproximateuniqueness}
\end{thm}
\en

\bn
\begin{thm}  Let $X$ be a compact, second countable, path connected,
metric space
and let $A$ be a unital
separable simple $\AH$-algebra with bounded dimension growth and real rank 
zero.
Suppose that $\alpha \in \KL (C(X), A)_{+, 1}$ 
and suppose that 
$\lambda : \TT (A) \rightarrow \TT (C(X))$ is an  
affine continuous map such that 
\begin{enumerate}
\item[(i)] the induced map 
$\Aff(\lambda) : \Aff(\TT (C(X)) \rightarrow \Aff(\TT (A))$ 
brings nonnegative nonzero 
functions to functions which are strictly positive at every point;
i.e.,  if 
\[
h \in C(X, \mathbb{R}) \cong \Aff(\TT (C(X))
\]
is  a nonnegative function
which is not the zero function then 
\[
\Aff(\lambda)(h)(\tau) > 0
\] 
for all
$\tau \in \TT (A)$; and       
\item[(ii)] for every $\tau \in \TT (A)$ and for every projection 
$p \in C(X) \otimes \K$,
\[
\lambda(\tau)(p) = \tau_{\ast}(\alpha([p])).
\]
\end{enumerate}
Then there exists a $*$-momomorphism $\phi : C(X) \rightarrow A$ such that\\
\[
\KL (\phi) = \alpha \makebox{ and } \TT(\phi) = \lambda.
\]
\label{thm:Existence}
\end{thm}

\begin{nproof}
Let $\{ \F_l \}_{l=1}^{\infty}$ be an increasing
sequence of finite subsets of 
the closed unit ball of $C(X)$ such that the union is dense in the closed
unit ball;  i.e., $\overline{C(X)_1} = 
\overline{ \bigcup_{l=1}^{\infty} \F_l }$.
Let $\{ \epsilon_l \}_{l=1}^{\infty}$ and 
$\{ \nu_l \}_{l=1}^{\infty}$ be strictly decreasing
sequences of strictly positive
real numbers 
such that $\sum_{l=1}^{\infty} \epsilon_l < \infty$ and 
$\sum_{l=1}^{\infty} \nu_l < \infty$, and such that 
$|f(x) - f(y) | < \epsilon_l/8$ for all $x,y \in X$ with 
$\dist (x,y) < \nu_l$ and for all $f \in \F_l$.  
For simplicity, we may
assume that $\nu_1 < 1$. 
For each $l$, 
\begin{enumerate}
\item take $s = l$;
\item take $\{ x_{l,1}, x_{l,2}, ..., x_{l, n_l} \}$ to 
be a $\nu_l/2$-dense
subset of $X$ for which $O_{l, j} \cap O_{l, k}
= \emptyset$ for $j \neq k$, where
$O_{l,j} = \{ x \in X : \dist (x, x_{l,j}) < \nu_l/(2l) \}$;  also, put
$\tilde{O}_{l,j} =_{df} \{ x \in X : \dist (x, x_{l, j}) < \nu_l/(4l) \} 
\subseteq
O_{l,j}$;   
\item take $\{ \sigma_l \}_{l=1}^{\infty}$ to be a strictly decreasing 
sequence of strictly positive real numbers such that 
$1/(2l) > \sigma_l > 0$ and for all $j$ for all $l$, 
$inf\{ \mu_{\lambda(\tau)}( \tilde{O}_{l ,j}) : 
\tau \in \TT ( A ) \} >
10 \sigma_l$.   
\end{enumerate}
For each $l$, put the above data 
into Theorem \ref{thm:Linapproximateuniqueness} 
and get the following:
\begin{enumerate}
\item[(a)]  a strictly decreasing sequence $\{ \gamma_l \}_{l=1}^{\infty}$ 
of strictly positive real numbers with 
$\sum_{l=1}^{\infty} \gamma_l < \infty$;
\item[(b)] an increasing sequence $\{ \G_l \}_{l=1}^{\infty}$ of finite
subsets of the closed unit ball of $C(X)$, which is dense in the closed
unit ball (i.e., $\overline{ \bigcup_{l=1}^{\infty} \G_l } = 
\overline{C(X)_1}$);  we may assume that $\F_l \subseteq \G_l$ for all $l$;
and
\item[(c)] a strictly decreasing sequence $\{ \delta_l \}_{l=1}^{\infty}$
of strictly positive real numbers such that $\sum_{l=1}^{\infty} \delta_l
< \infty$. 
\item[(d)] an increasing sequence $\{ \P_l \}_{l=1}^{\infty}$ of 
finite subsets of $\mathbb{P}(C(X))$.
\item[(e)]  for each $l,j$, let $f_{l,j}$ be a real-valued function such
that i. $0 \leq f_{l,j} \leq 1$, ii. $f_{l,j} = 1$ on $\tilde{O}_{l,j}$
and iii. the support of $f_{l,j}$ is contained in $O_{l,j}$; 
then, expanding $\G_l$ if necessary, we may assume that 
$\{ f_{l,j} : 1 \leq j \leq n_l \} \subseteq \G_l$; 
\item[(f)] further expanding the $\G_l$s  and contracting the $\gamma_l$s and 
$\delta_l$s if necessary, we may assume that for each $l$ and for every
unital $C^*$-algebra $\mathcal{C}$, if $\kappa_1, \kappa_2 : C(X) 
\rightarrow  \mathcal{C}$ are $\G_l - \delta_l$-multiplicative 
c.p.c. maps, then they are defined on $\P_l$;  and if, in 
addition, $\kappa_1 (f)$ is within $\gamma_l$ of $\kappa_2(f)$ for all
$f \in \G_l$ then
\[
[ \kappa_1 ] | \P_l = [ \kappa_2 ] | \P_l.
\] 
\end{enumerate}

     Hence, for each $l$, for the given data $\F_l, \epsilon_l, \nu_l,
s=l, \{ x_{l,1}, x_{l,2}, ..., x_{l, n_l} \}, \sigma_l$, the quantities
$\gamma_l, \G_l, \P_l$ satisfy the conclusion of Theorem 
\ref{thm:Linapproximateuniqueness}.

     Since $A$ is a unital simple $\AH$-algebra with bounded dimension growth
and real rank zero, 
$A$ is (the norm closure of)
an increasing union $A = \overline{\bigcup_{m=1}^{\infty} A_m}$, where
each $A_m$ is a direct sum of
unital homogeneous $C^*$-algebras with spectra being 
finite CW complexes with topological dimension less than or equal to three,
 and
where (we may assume that) the inclusions are unital and injective.

     Hence, we have that $\TT ( A )$ is the inverse limit of tracial state 
spaces:
\[
\TT (A_1) \leftarrow \TT (A_2) \leftarrow \TT (A_3) \leftarrow ... \leftarrow 
\TT ( A ).
\]
Hence, $\Aff (\TT ( A ))$ is the direct limit of complete order unit spaces
$\Aff (\TT (A_m))$:
\[
\Aff (\TT (A_1)) \rightarrow \Aff (\TT (A_2)) \rightarrow 
\Aff (\TT (A_3)) \rightarrow ...
\rightarrow \Aff (\TT (  A )).
\]
Note that if the spectrum of $A_m$ is $Y$ then $\Aff (\TT (A_m))$ is isomorphic,
as an order unit space, to $C(Y, \mathbb{R})$ (the real-valued continuous
functions on $Y$). 
Note also that $\lambda$ induces a morphism (of order unit spaces)
$\Aff (\lambda) : \Aff (\TT (C(X)))   \rightarrow \Aff (\TT ( A ))$;  and 
by the argument of \cite{Fields} Theorem 25.1.1, 
for every $\epsilon* >0 $ and for every finite subset $\F*$ of $\Aff 
(\TT (C(X)))$,
there exists an integer $N* \geq 1$ such that for every $n \geq N*$, there 
is a morphism (of order unit spaces) $\Aff (\lambda*) : \Aff (\TT (C(X))) 
\rightarrow
\Aff (\TT (A_n))$ where $\Aff (\lambda)(f)$ is within $\epsilon*$ of 
$\Aff (\lambda*)(f)$ for every $f \in \F*$ (viewing $\Aff (\lambda*)$ as a map
with codomain $\Aff (\TT ( A ))$ by taking the natural composition etc.).    

    For each $l$, let $\G_l'$ consist of positive elements of norm less 
than or equal to one such that $\G_l' =_{df} \{ a_1, a_2, a_3, a_4 \geq 0 : 
f = a_1 - a_2 + i (a_3 - a_4), f \in \G_l
 \makebox{  and  } a_1 a_2 = a_3 a_4 = 0 \}$. 
Hence, $\{ \G_l' \}_{l=1}^{\infty}$ is an increasing sequence of finite
subsets of $C( X, \mathbb{R})$ such that the union is dense in the 
closed unit ball of the positive elements in $C(X, \mathbb{R})$.  

     We construct sequences $\{ \psi_l \}_{l=1}^{\infty}$
and a subsequence $\{ N_l \}_{l=1}^{\infty}$ of the integers such that 
\begin{enumerate}
\item  $\psi_l : C(X_1) \rightarrow A_{N_l}$ is a unital 
c.p.c. $\G_l - \delta_l$ - multiplicative map; 
\item $\tau \circ \psi_l (f)$ is within $\gamma_l/4$ of $\lambda(\tau)(f)$
for all $f \in \G_l$;
\item   $\mu_{\tau \circ \psi_l}( O_{l',j}) > \sigma_{l'}\nu_{l'}$ for
all $l' \leq l$, for all $j$ and all $\tau \in \TT ( A )$;  and  
\item $\alpha | \P_l = [\psi_l ] | \P_l$.  
\end{enumerate}
Let us collectively denote the above conditions by ``(*)".

To simplify notation, let us denote $\rho_l =_{df} \mmin \{ \gamma_l /100, 
\sigma_l/100 \}$.  
Fix $l$.  
We now construct $\psi_l$.    
By the same argument as that of \cite{Fields} Theorem 25.1.1,
there exists an integer $M_1$ such that for all $n \geq M_1$ there is a 
morphism (of complete order unit spaces) $\Aff (\lambda_n*) : \Aff (\TT (C(X)))
\rightarrow \Aff (\TT (A_n))$ such that
$\Aff (\lambda)(f)$ is within 
$\rho_l$ of $\Aff (\lambda_n*)(f)$ for all 
$f \in \G_l'$ (of course, we mean the appropriate composition of 
$\Aff (\lambda_n*)$ with etc. to get a map with codomain $\Aff (\TT (A))$).
Since $ A $ is a simple unital $\AH$-algebra with bounded dimension
growth and real rank zero,  and 
by Theorem 
\ref{thm:Linexistence}, let $\phi_l : C(X) \rightarrow A$ be a unital,
c.p.c., almost multiplicative map and  
let $q$ be a projection in $ A $ such that 
\begin{enumerate}
\item[i.]  there is an integer $M_2 \geq M_1$ such that  
$q \in A_{M_2}$ and $\tau'(q) <  \rho_l$
for every $\tau' \in \TT (A_{M_2})$;  
\item[ii.] the map $C(X) \rightarrow A : f \mapsto 
q \phi_l (f) q$ is a c.p.c. 
$\G_l - \delta_l$-multiplicative map;
\item[iii.]  there is a 
finite dimensional $C^*$-subalgebra $F \subseteq A_{M_2}$ with
$\be_F = 1_A - q$;  
\item[iv.]  there is a unital finite-dimensional $*$-homomorphism
$h : C(X) \rightarrow F$ such that 
$\phi_l(f)$ is within $\rho_l$ of
$q \phi_l (f) q + h(f)$ for all $f \in \G_1$;  
\item[v.] $\alpha | \P_l = [ \phi_l ] | \P_l$.
\end{enumerate}

     Let $\phi_{l,1} : C(X) \rightarrow A$ be 
the unital 
c.p.c. $\G_l - \delta_l$-multiplicative map given by 
$\phi_{l,1}(f) =_{df} q \phi_l(f) q + h(f)$ for all $f \in C(X)$. 
Note that it follows, from the above conditions and by the definition
of the $\G_l$ (and the remarks surrounding it) that    
$\phi_{l,1}$ is well-defined on $\P_l$ and
\[
\alpha| \P_l = [\phi_l]| \P_l = [ \phi_{l,1}] | \P_l.
\]
 
     Put $\G_l'$ and $\rho_l$ into Theorem \ref{thm:Litraces} 
to get the integer $L$ (which is the $N$ in the statement
of Theorem \ref{thm:Litraces}). 
Note that $A_{M_2}$ has the form 
$A_{M_2} = r_1 \M_{k_1} (C(Y_1)) r_1 \oplus r_2 \M_{k_2} (C(Y_2)) r_2 \oplus
... \oplus r_a \M_{k_a} (C(Y_a)) r_a$ where each $Y_j$ is a connected 
finite CW-complex with dimension less than or equal to three.
Suppose that $F = F_1 \oplus F_2 \oplus ... \oplus F_a$ where for every $j$,
$F_j$ is (a finite dimensional $*$-algebra) contained in 
$r_j \M_{k_j} (C(Y_j)) r_j$.   
Since $A$ is a simple unital $\AH$-algebra with bounded dimension
growth and real rank zero, 
moving up building blocks if necessary (see \cite{EG}), 
we may assume
that for each $j$,  there is a trivial 
projection $t_j \in r_j \M_{k_j} (C(Y_j)) r_j$  (i.e., a projection 
corresponding to a trivial vector bundle over $Y_j$) 
such that (a) 
$L t_j$ is Murray-von Neumann equivalent to a subprojection of
$\be_{F_j}$ in $r_j \M_{k_j} (C(Y_j)) r_j$ and 
(b) $\tau'( [\be_{F_j} - L t_j ]) < \rho_l$ 
for all $\tau' \in \TT (F_j)$.

     Now recall that 
since $M_2 \geq M_1$, the morphism $\Aff (\lambda_{M_2}*) : 
\Aff (\TT (C(X)) \rightarrow \Aff (\TT (A_{M_2}))$ (of order unit spaces) 
is such 
that $\Aff (\lambda)(f)$ is within $\rho_l$
of $\Aff (\lambda_{M_2}*)(f)$ for all $f \in \G_l'$.
Note that $\Aff (\TT (C(X))) \cong C(X, \mathbb{R})$ and 
$$\Aff (\TT (A_{M_2})) \cong C(Y_1, \R) \oplus C(Y_2, \R) \oplus .... 
\oplus C(Y_a, \R),$$ 
where the isomorphisms are isomorphisms between order unit spaces.
For each $j$, there is the natural projection map
$\pi_j : C(Y_1, \R) \oplus C(Y_2, \R) \oplus ... \oplus C(Y_j, \R) 
\rightarrow C(Y_j, \R)$
which is a (surjective) morphism of order unit spaces.
Hence, for each $j$, we get a morphism (or, unital positive map)
$\xi_j =_{df} \pi_j \circ \Aff (\lambda_{M_2}*) : C(X) \rightarrow C(Y_j)$.
By the definition of $L$ and by Theorem \ref{thm:Litraces},
let $\chi_{j,1}, \chi_{j,2}, ..., \chi_{j,L} : 
C(X) \rightarrow C(Y_j)$ be unital
$*$-homomorphisms such that
each $\chi_{j,k}$ is homotopy equivalent to a 
point-evaluation $*$-homomorphism
and such that 
\[
\| \xi_j(f) - (1/L) \sum_{k=1}^L \Aff (\TT (\chi_{j,k})) (f) \| < 
\rho_l 
\]
for all $f \in \G_l'$. 

     For each $j$, 
let $t_{j,1}, t_{j,2}, ..., t_{j,L}$ be the $L$ pairwise orthogonal
subprojections of 
$\be_{F_j}$, each of which is Murray-von Neumann equivalent (in $F_j$) to
$t_j$; and  
let $h_j' : C(X) \rightarrow \be_{F_j} \M_{k_j}(C(Y_j)) \be_{F_j}$ 
be the unital 
$*$-homomorphism given by
$h_j' : f \mapsto \sum_{k=1}^L \chi_{j,k}(f) \otimes t_{j,k} + 
f(z_j)(\be_{F_j} - \sum_{k=1}^L t_{j,k})$, for
some point $z_j \in X$. 
Let $h' : C(X) \rightarrow \be_F \A_{M_2} \be_F$ be given
by $h' =_{df} h_1' \oplus h_2' \oplus ... \oplus h_a'$.
If $h = h_1 \oplus h_2 \oplus ... \oplus h_a$ where each $h_j : C(X) 
\rightarrow F_j$ is a unital $*$-homomorphism then
$h_j$ is homotopy equivalent to $h_j'$, for every $j$. 
Hence, $h$ is homotopy equivalent to $h'$.
Let $\psi_l : C(X) \rightarrow A$ be given by
$\psi_l : f \mapsto q \phi_l(f) q  + h'(f)$ for all $f \in C(X)$.
Note that $\psi_l$ is a unital c.p.c.
$\G_l - \delta_l$-multiplicative map.
Then, by the definition of $\G_l$ (and the remarks surrounding it),
$\psi_l$ is well-defined on $\P_l$ and 
\[
[\psi_l] | \P_l = [ \phi_{l,1} ] | \P_l.
\]
Hence,
\[
[\psi_l] | \P_l = [ \phi_l ] | \P_l = \alpha| \P_l.
\]

    Fix $1 \leq j \leq a$.  For $\tau' \in \TT (r_j \M_{k_j} ( C(Y_j)) r_j)$ and
$f \in \G_l'$, 
\begin{eqnarray*} 
\lefteqn{ | ((1/L) \sum_{k=1}^L \tau'(\chi_{j,k}(f)) - 
\tau'(h_j'(f)) |}   \\ 
& = &  | ((1/L) \sum_{k=1}^L \tau'(\chi_{j,k}(f)) -  
(\sum_{k=1}^L \tau'(\chi_{j,k}(f) \otimes t_{j,k}) + 
\tau'(f(z_j)(\be_{F_j} - 
\sum_{k=1}^L t_{j,k})) | \\ 
& < &    | (1/L)   \sum_{k=1}^L \tau'(\chi_{j,k}(f)) - 
\sum_{k=1}^L \tau'( \chi_{j,k}(f) \otimes
t_{j,k}) | + \rho_l  \\
& < & 3 \rho_l.
\end{eqnarray*}
Hence, for $f \in \G_l'$, and for 
$\tau' \in \TT (r_j \M_{k_j}(C(Y_j)) r_j)$,  
\begin{eqnarray*}
\lefteqn{| \tau'(\xi_j (f))  - \tau'(h_j'(f)) |}  \\ 
& = & | \tau'(\xi_j(f)) - (1/L) \sum_{k=1}^L \tau'(\chi_{j,k}(f)) | + 
| (1/L) \sum_{k=1}^L \tau'(\chi_{j,k}(f)) - \tau'(h_j'(f)) |  \\ 
& < & \rho_l + | (1/L) \sum_{k=1}^L (\tau'(\chi_{j,k}(f)) - 
\tau'(h_j'(f))) |  \\
& < & 4 \rho_l. 
\end{eqnarray*}

    Hence, for $f \in \G_l'$ and for $\tau \in \TT ( A )$,
\begin{eqnarray*}
\lefteqn{| \lambda(\tau)(f) - \tau(\psi_l(f)) |} \\
& \leq & | \lambda(\tau)(f) - \tau(\Aff (\lambda_{M_2}*)(f)) | + 
| \tau(\Aff (\lambda_{M_2}*)(f)) - \tau(\psi_l(f)) | \\ 
& < & \rho_l + | \tau(\Aff (\lambda_{M_2}*)(f)) - \tau(\psi_l(f)) |  \\ 
& = & \rho_l + | \tau(\Aff (\lambda_{M_2}*)(f)) - (\tau(q \phi_l(f) q) + 
\tau(h'(f)) |  \\ 
& < & 2\rho_l + | \tau(\Aff (\lambda_{M_2}*)(f)) - \tau(h'(f)) |  \\   
& < & 6 \rho_l.
\end{eqnarray*}

     Since for $l_0 \leq l$,
$\{ f_{l_0,j} : 1 \leq j \leq n_{l_0} \} \subseteq \G_l'$ and
since $\sigma_l < \sigma_{l_0}$, it 
follows that $\tau(\psi_l (f_{l_0,j})) > \sigma_{l_0}$ 
for all $\tau \in \TT ( A  )$ and for all $j$. 
Hence, $\mu_{\tau \circ \psi_l}(O_{l_0,j}) > \sigma_{l_0} \nu_{1_0}$ 
for all $l_0 \leq l$ and all $j$.     

Also, from the above and definitions of $\G_l$ and $\G_l'$, it follows that 
for all $\tau \in \TT (A )$ and for all $f \in \G_l$,
$| \lambda(\tau)(f) - \tau(\psi_l (f)) | < 24 \rho_l < \gamma_l/4$.   

  Hence, we have constructed a sequence $\{ \psi_l \}_{l=1}^{\infty}$ that 
satisfies the conditions in (*). 
Hence, by (*) and Theorem \ref{thm:Linapproximateuniqueness},  
let $u_1$ be a unitary in $A$ such that 
$u_1 \psi_2(f) u_1^*$ is within $\epsilon_1$ of $\psi_1(f)$ for all
$f \in \F_1$.
Again by (*) and Theorem  \ref{thm:Linapproximateuniqueness},  
let $u_2$ be a unitary in $A$ such that
$u_2 \psi_3(f) u_2^*$ is within $\epsilon_2$ of $u_1 \psi_2(f) u_1^*$
for all $f \in \F_2$.  Repeating this process, we get a sequence
$\{ u_l \}_{l=0}^{\infty}$ of unitaries in $A$ (taking $u_0 = 1_A$)
such that for every $l$,
$u_l \psi_{l+1}(f) u_l^*$ is within $\epsilon_l$ of 
$u_{l-1} \psi_l (f) u_{l-1}^*$ for all $f \in \F_l$.
Since the union of the $\F_l$s is dense in the closed unit ball of
$C(X)$ and since $\sum_{l=1}^{\infty} \epsilon_l < \infty$,
the sequence $\{ u_l \psi_{l+1} u_l^* \}_{l=1}^{\infty}$ must converge
pointwise to a unital 
$*$-homomorphism $\Phi : C(X) \rightarrow  A $.
And by our construction, $\lambda = \TT (\Phi)$.
Moreover, $\KL (\Phi) = \KL (\alpha)$.  
\end{nproof}
\en

\bn
    From Theorem \ref{thm:Existence}
and \cite{Linembedding} Theorem 4.5, we get Theorem \ref{thm:MainExistence}. 


Finally, since a finite
CW complex is the finite union of pairwise disjoint,
path connected, second countable, compact metric spaces, 
we can replace $X$ in Theorem~\ref{thm:MainExistence} by an arbitrary finite
CW complex.

 
\en 

\section{Tracially $\AF$ embeddings}

\bn
\label{commutative-tr0-embedding}
It is the aim of this section to prove the lemma below. Our argument is inspired by the methods developed in \cite{Wi5} and \cite{Wi6}.  

\begin{lms}
Let $Y$ be a compact metrizable space and $A$ a simple, separable, unital, $\Zh$-stable $C^{*}$-algebra with real rank zero. Suppose $\theta: \Ch(Y) \to A$ is a unital $*$-homomorphism; let a finite subset $\Fh \subset \Ch(Y)$ and $\varepsilon>0$ be given.\\
Then, there is a commutative finite-dimensional $C^{*}$-subalgebra $B \subset A$ such that
\begin{itemize}
\item[(i)] $\|[\theta(a),\be_{B}]\|<\varepsilon \; \forall \, a \in \Fh$
\item[(ii)] $\dist(\be_{B}\theta(a)\be_{B},B)< \varepsilon \; \forall \, a \in \Fh$
\item[(iii)] $\tau(\be_{B}) > 1 - \varepsilon \; \forall \, \tau \in \TT(A)$.
\end{itemize}
\end{lms}

A $*$-monomorphism that satisfies the conclusion of the above lemma will
be called ``tracially $\AF$-embeddings" or ``$\TAF$-embeddings" (see Definition~\ref{d-tr0-embedding} below).
\en

\bn
\begin{nots}
For a $C^{*}$-algebra $A$, we denote by $A_{\infty}$ the quotient $\prod_{\N} A / \bigoplus_{\N} A$.
\end{nots}
\en

\bn
\label{commutative-AF-approximation}
\begin{props}
Let $Y$ be a compact metrizable space. Then, there are a zero-dimensional compact metrizable space $X$, a unital embedding $\nu: \Ch(Y) \to \Ch(X)$ and a commutative diagram as follows:
\begin{equation}
\label{small-diagram}
\xymatrix{
& \Ch(Y) \ar[ddl]_{(\psi_{l})_{\N}} \ar[d]^{\nu} \\
& \Ch(X) \ar[d]^{\kappa} \\
\prod_{\N} \C^{k_{l}} \ar[r]^-{q} & \prod_{\N} \C^{k_{l}} / \bigoplus_{\N} \C^{k_{l}}
}
\end{equation}
Here, $\Ch(X) = \lim_{\to} (\C^{k_{l}}, \beta_{l})$ is a representation of $\Ch(X)$ as a unital $\AF$ algebra, $\kappa: \Ch(X) \to \prod \C^{k_{l}}/\bigoplus \C^{k_{l}}$ the canonical inclusion of the inductive limit, $q$ the  quotient map and $(\psi_{l})_{l \in \N}$ a sequence of unital $*$-homomorphisms. \\
Moreover, there is a sequence of c.p.c.\ maps 
\[
\varphi_{l}: \C^{k_{l}} \to \Ch(Y)
\]
such that $q \circ ((\varphi_{l} \psi_{l})_{l \in \N}) = \iota_{\Ch(Y)}$, where $\iota_{\Ch(Y)}: \Ch(Y) \to \prod_{\N} \Ch(Y)$ is the canonical embedding (as constant sequences) and $q$ also denotes the quotient map $\prod_{\N} \Ch(Y) \to \Ch(Y)_{\infty}$. If $\dim Y= n<\infty$, we may choose the $\varphi_{l}$ to be $n$-decomposable in the sense of \cite[Definition~2.2]{KW}.
\end{props}

\begin{nproof}
Choose a sequence $(\Uh_{l}=\{U_{l,1}, \ldots,U_{l,k_{l}}\})_{l \in \N}$ of open coverings  of $Y$ with the following properties:
\begin{enumerate}
\item[(i)] $\Uh_{l+1}$ refines $\Uh_{l}$ for all $l$
\item[(ii)] each $U_{l,j}$ contains some element $y_{l,j}$ of $Y$
\item[(iii)] for some fixed metric on $Y$, $\max_{j} \{\diam(U_{l,j}) \} \stackrel{l \to \infty}{\longrightarrow}0$.
\end{enumerate}
Define $*$-homomorphisms $\psi_{l}: \Ch(Y) \to \C^{k_{l}}$ by $\psi_{l}:= \bigoplus_{j=1}^{k_{l}} \ev_{y_{l,j}}$. Choose partitions of unity subordinate to $\Uh_{l}$; interpret these as c.p.c.\ maps $\varphi_{l}:\C^{k_{l}} \to \Ch(Y)$. It is clear from (iii) that $\varphi_{l}\psi_{l} \to \id_{\Ch(Y)}$ pointwise, whence $q((\varphi_{l}\psi_{l})_{l\in \N}) = \iota_{\Ch(Y)}$. Note also that, if $\dim Y= n <\infty$, then the $\Uh_{l}$ -- and therefore also the $\varphi_{l}$ -- may be chosen to be $n$-decomposable.\\
Since $\Uh_{l+1}$ refines $\Uh_{l}$ for each $l$, we may choose unital $*$-homomorphisms $\beta_{l}: \C^{k_{l}} \to \C^{k_{l+1}}$ such that, if the $i$-th component of $\beta_{l}(e_{l,j})$ is nonzero, then $U_{l+1,i} \subset U_{l,j}$ (where $e_{l,j}$ denotes the $j$-th canonical generator of $\C^{k_{l}}$). In other words, for each pair $(l+1,i)$ we choose a pair $(l,j)$ such that  $U_{l+1,i} \subset U_{l,j}$ and let this assignment represent an arrow in the Bratteli diagram $\C^{k_{l}} \to \C^{k_{l+1}}$. The inductive limit $\lim_{\to}(\C^{k_{l}},\beta_{l})$ is a commutative unital $\AF$ algebra, hence of the form $\Ch(X)$ for some compact zero-dimensional space $X$; let $\kappa: \Ch(X) \to \prod \C^{k_{l}}/\bigoplus \C^{k_{l}}$ denote the canonical embedding. \\
For each $l \in \N$ we have a $*$-homomorphism 
\[
\nu_{l}: \Ch(Y) \stackrel{\psi_{l}}{\longrightarrow} \C^{k_{l}} \stackrel{\beta_{l,\infty}}{\longrightarrow} \Ch(X) \, ;
\]
using (iii) it is now straightforward to check that this sequence of maps is approximately multiplicative in the sense of \cite{BK1} and induces a $*$-homomorphism $\nu : \Ch(Y) \to \Ch(X) $  which makes our diagram (\ref{small-diagram}) commute. 
\end{nproof}
\en

\bn
\label{order-0-excision}
The next lemma uses some technical results of \cite{Wi5}; it is 
essentially contained in the proof of Theorem 4.1 of that paper. 

\begin{lms}
Let $A$ and $B$ be separable $C^{*}$-algebras, $A$ unital and $\Zh$-stable with real rank zero. Let $n,k \in \N$ and  $p \in A$ be a projection. Consider c.p.c.\ maps
\[
B \stackrel{\psi}{\longrightarrow} \C^{k} \stackrel{\varphi}{\longrightarrow} pAp
\]
such that $\varphi$ is $n$-decomposable for some $n \in \N$ and satisfies $\|p - \varphi(\be_{\C^{k}})\|\le \eta$ for some $\eta > 0$. \\
Then, there is a sequence of $*$-homomorphisms $\varrho_{m}:\C^{k} \to pAp$, $m \in \N$, such that
\[
\limsup_{m \to \infty} \|[\varrho_{m}(\be_{\C^{k}}), \varphi \psi(b)]\| \le 10(n+1) \cdot \|\varphi \psi(b^{2}) - \varphi \psi(b)^{2}\|^{\halb} \; \forall \, b \in B_{+} \, ,
\]
\[
\limsup_{m \to \infty} \|\varrho_{m}(\be_{\C^{k}}) \varphi \psi(b) - \varrho_{m} \psi(b)\| \le 5(n+1) \cdot \|\varphi \psi(b^{2}) - \varphi \psi(b)^{2}\|^{\halb} \; \forall \, b \in B_{+} 
\]
and
\[
\tau(\varrho_{m}(\be_{\C^{k}})) >  \left(\frac{1}{5(n+1)} - \eta \right) \cdot \tau(p)  \; \forall \, \tau \in \TT (A), \, m \in \N \, .
\]
\end{lms}
\en

\bn
\label{expansions}
Let $(j \mapsto l_{j})$ be a surjective and decreasing map $\N \to \N$. We call the sequence $(l_{j})_{j \in \N}$ an expansion of the sequence $(l)_{l \in \N}$.  For a sequence $(M_{l})_{l \in \N}$ of maps, algebras, etc., we call the sequence $(M_{l_{j}})_{j \in \N}$ an expansion of $(M_{l})_{l \in \N}$. Similarly, from an inductive system $(A_{l},\alpha_{l})_{l \in \N}$ of $C^{*}$-algebras we obtain an expansion $(A_{l_{j}}, \bar{\alpha}_{j})_{j \in \N}$, where $\bar{\alpha}_{j}=\alpha_{l_{j}}$ if $l_{j} \neq l_{j+1}$, and $\bar{\alpha}_{j} = \id_{A_{l_{j}}}$ if $l_{j} = l_{j+1}$. It is clear that $\lim_{l \to \infty} A_{l} \cong \lim_{j \to \infty}A_{l_{j}}$. 
\en

\bn
\label{theta-approximation}
\begin{props}
Let $A$ be a simple, separable, unital, $\Zh$-stable $C^{*}$-algebra with real rank zero; let $Y$ be a compact metrizable space with $\dim Y = n < \infty$. Suppose $\theta: \Ch(Y) \to A_{\infty}$ is a $*$-homomorphism and $(p_{l})_{l\in \N} \in \prod A$ is a  sequence of projections lifting $p:= \theta(\be_{Y})$. \\
Then, there is a commutative diagram as follows:
\begin{equation}
\xymatrix{
& \Ch(Y) \ar[d]_{\nu} \ar[ddl]_{(\psi_{l})_{\N}} \ar[dddr]^{h \theta(\,.\,)h}& \\
& \Ch(X) \ar[d]_{\kappa} & \\
\prod_{\N} \C^{k_{l}} \ar[r]^-{q} \ar[d]_{(\varrho_{l})_{\N}} & \prod_{\N} \C^{k_{l}}/\bigoplus_{\N} \C^{k_{l}} \ar[d]_{\varrho} & \\
\prod_{\N} p_{l}Ap_{l} \ar[r]^{q} & pA_{\infty} p \ar[r]^{\subset} & A_{\infty} 
}
\label{large-diagram}
\end{equation}
Here, the upper left triangle has the properties of diagram (\ref{small-diagram}) in Proposition \ref{commutative-AF-approximation} (in fact, it is an expansion of (\ref{small-diagram})),  the $\varrho_{l}$ and $\varrho$ are $*$-homomorphisms, 
\[
h:= \varrho \kappa \nu(\be_{Y}) \in \theta(\Ch(Y))' \cap A_{\infty}
\]
and the map $h \theta(\, . \, )h$ is a $*$-homomorphism. Moreover, the $\varrho_{l}$ satisfy 
\begin{equation}
\label{theta-approximation-estimate}
\liminf_{l \to \infty} \left( \tau(\varrho_{l} \psi_{l}(\be_{Y})) - \frac{1}{5(n+1)} \cdot \tau(p_{l}) \right) \ge 0 \; \forall  \, \tau \in \TT (A) \, .
\end{equation}
\end{props}

\begin{nproof} 
Given $Y$, apply Proposition \ref{commutative-AF-approximation} to obtain a diagram 
\begin{equation*}
\xymatrix{
& \Ch(Y) \ar[ddl]_{(\tilde{\psi}_{l})_{\N}} \ar[d]^{\tilde{\nu}} \\
& \Ch(X) \ar[d]^{\tilde{\kappa}} \\
\prod \C^{k_{l}} \ar[r]^-{q} & \prod \C^{k_{l}} / \bigoplus \C^{k_{l}}
}
\end{equation*}
and $n$-decomposable maps $\tilde{\varphi}_{l}: \C^{k_{l}} \to \Ch(Y)$. Set 
\[
\bar{\varphi}_{l} := \theta \circ \tilde{\varphi}_{l} \; \forall \, l \in \N \, ;
\] 
these are $n$-decomposable c.p.c.\ maps $\C^{k_{l}} \to p A_{\infty} p$ satisfying
\[
\bar{\varphi}_{l}\circ \tilde{\psi}_{l} (b) \stackrel{l \to \infty}{\longrightarrow} \theta(b) \; \forall \, b \in \Ch(Y) \, .
\]
By \cite{KW}, Remark 2.4 (cf.\ also \cite{Wi2}, Proposition 1.2.4), each $\bar{\varphi}_{l}$ may be lifted to a sequence of $n$-decomposable c.p.c.\ maps 
\[
\bar{\varphi}_{l,j} : \C^{k_{l}} \to p_{j}Ap_{j}, \; j \in \N \, .
\]
A diagonal sequence argument now yields an expansion $(l_{j})_{j \in \N}$ of the sequence $(l)_{l \in \N}$ such that 
\[
\varphi_{j} := \bar{\varphi}_{l_{j},j} : \C^{k_{l_{j}}} \to p_{j}Ap_{j}
\]
and
\[
\psi_{j}:= \tilde{\psi}_{l_{j}}
\]
satisfying
\begin{equation}
\label{theta-factorization}
q \circ ((\varphi_{j} \circ \psi_{j})_{j \in \N}) = \theta \, .
\end{equation}
Note that the inductive system $\Ch(X) \cong \lim_{\to} \C^{k_{l}}$ and its expansion $\lim_{\to} \C^{k_{l_{j}}}$ are isomorphic by \ref{expansions}. Let $\nu$ denote the composition of this isomorphism with $\tilde{\nu}$ and write $\kappa$ for the natural inclusion of $\lim_{\to} \C^{k_{l_{j}}}$ into $\prod \C^{k_{l_{j}}}/\bigoplus \C^{k_{l_{j}}}$. This yields a commutative diagram
\begin{equation*}
\xymatrix{
& \Ch(Y) \ar[ddl]_{(\psi_{j})_{\N}} \ar[d]^{\nu} \\
& \Ch(X) \ar[d]^{\kappa} \\
\prod \C^{k_{l_{j}}} \ar[r]^-{q} & \prod \C^{k_{l_{j}}} / \bigoplus \C^{k_{l_{j}}} 
}
\end{equation*}
with the same properties as (\ref{small-diagram}). Misusing notation, from now on we will write $\C^{k_{j}}$ in place of $\C^{k_{l_{j}}}$; this should not cause confusion.\\
For each $j \in \N$, we may apply Lemma \ref{order-0-excision} with $\Ch(Y)$ in place of $B$ and $p_{j},\C^{k_{j}}, \psi_{j},\varphi_{j}$ and $\eta_{j}:= \|p_{j}- \varphi(\be_{\C^{k_{j}}})\|$ in place of $p,\C^{k},\psi,\varphi$ and $\eta$, respectively, to obtain a sequence of $*$-homomorphisms 
\[
\varrho_{j,m}:\C^{k_{j}} \to p_{j}Ap_{j} 
\]
satisfying the assertions of \ref{order-0-excision}. Note that 
\[
\|p_{j} - \varphi_{j}(\be_{\C^{k_{j}}})\| \stackrel{j \to \infty}{\longrightarrow} 0
\]
and
\[
\|\varphi_{j} \psi_{j}(b^{2}) - \varphi_{j} \psi_{j}(b)^{2}\| \stackrel{j \to \infty}{\longrightarrow} \; \forall \, b \in \Ch(Y) \, .
\]
Therefore, again by a diagonal sequence argument (and separability of $\Ch(Y)$) we may choose a sequence $(m_{j})_{j \in \N} \subset \N$ such that for the $*$-homomorphisms 
\[
\varrho_{j}:= \varrho_{j,m_{j}} : \C^{k_{j}} \to p_{j}Ap_{j} , \, j \in \N
\]
we have 
\begin{equation}
\label{central-rho}
\|[\varrho_{j}(\be_{\C^{k_{j}}}), \varphi_{j}\psi_{j}(b)]\| \stackrel{j \to \infty}{\longrightarrow} 0 \; \forall \, b \in \Ch(Y)\, ,
\end{equation}
\begin{equation}
\label{excising-rho}
\|\varrho_{j}(\be_{\C^{k_{j}}}) \varphi_{j}\psi_{j}(b) - \varrho_{j}\psi(b)\| \stackrel{j \to \infty}{\longrightarrow} 0 \; \forall \, b \in \Ch(Y)
\end{equation}
and
\begin{eqnarray}
\label{tau-rho}
\lefteqn{ \liminf_{j} \left(\tau(\varrho_{j}\psi_{j}(\be_{Y})) - \frac{1}{5(n+1)} \cdot \tau(p_{j}) \right)} \nonumber \\
&=& \liminf_{j} \left(\tau(\varrho_{j}(\be_{\C^{k_{j}}})) - \frac{1}{5(n+1)} \cdot \tau(p_{j}) \right) \nonumber \\
& \ge &   0 \; \forall \, \tau \in \TT (A) \, .
\end{eqnarray}
The $\varrho_{j}$ now yield a commutative diagram 
\begin{equation*}
\xymatrix{
\prod \C^{k_{j}} \ar[r]^-{q} \ar[d]^{(\varrho_{j})_{\N}} & \prod \C^{k_{j}}/\bigoplus \C^{k_{j}} \ar[d]^{\varrho} \\
\prod p_{j}Ap_{j} \ar[r]^-{q} & p A_{\infty} p \, ,
}
\end{equation*}
where $\varrho$ denotes the map induced by the $\varrho_{j}$. The above diagram will be the lower left rectangle of (\ref{large-diagram}). It follows from (\ref{theta-factorization})  and 
(\ref{central-rho}) that 
\[
[\varrho(\be_{\prod \C^{k_{j}}/\bigoplus \C^{k_{j}}}), \theta(b)] = 0 \; \forall \, b \in \Ch(Y) \, .
\]
Since the $\psi_{j}$ are unital, this implies that 
\[
h = \varrho  \kappa \nu(\be_{Y}) = q((\varrho_{j}\psi_{j}(\be_{Y}))_{j \in \N})  = \varrho(\be_{\prod \C^{k_{j}}/\bigoplus \C^{k_{j}}}) \in \theta(\Ch(Y))' \cap A_{\infty} \, .
\]
Note that (\ref{theta-factorization}), (\ref{central-rho}) and (\ref{excising-rho}) also imply that
\[
h \theta(b) h = \varrho  \kappa  \nu (\be_{Y}) \theta(b)= q((\varrho_{j}\psi_{j}(b))_{j \in \N}) \; \forall \, b \in \Ch(Y) \, .
\]
Exchanging the indices $j$ for $l$, we have constructed all the ingredients for the diagram (\ref{large-diagram}); above we have proved commutativity of this diagram. The proposition's statement about traces is just (\ref{tau-rho}).
\end{nproof}
\en

\bn
\label{cantor-interpolation}
\begin{props}
Let $Y$ be a compact metrizable space with $\dim Y=n < \infty$. Let $A$ be a simple, separable, unital, $\Zh$-stable $C^{*}$-algebra with real rank zero. Suppose $\theta:\Ch(Y) \to A$ is a unital $*$-homomorphism and let $\varepsilon>0$ be given.\\
Then, there are a commutative $\AF$ algebra $\Ch(X) \subset A_{\infty}$,  a sequence $(p_{l})_{l \in \N} \in \prod_{\N} A$ of projections lifting $p:= \be_{X} \in A_{\infty}$ and $l_{0} \in \N$ such that 
\begin{itemize}
\item[(i)] $p \in \theta(\Ch(Y))' \cap A_{\infty}$ 
\item[(ii)] $p \theta(\Ch(Y)) \subset \Ch(X)$
\item[(iii)] $\tau(p_{l})> 1- \varepsilon \; \forall \, l \ge l_{0}, \, \tau \in \TT (A)$.
\end{itemize}
\end{props}

\begin{nproof}
For convenience, we set 
\[
\mu:= \frac{1}{5(n+1)} \, ;
\]
let $\theta^{(0)}:= \theta$ and $p_{l}^{(0)}:=\be_{A}$ for $l \in \N$ and $p^{(0)}:= \be_{A_{\infty}}$. Provided that, for some $m \in \N$, we have constructed a unital $*$-homomorphism $\theta^{(m)}: \Ch(Y) \to p^{(m)}A_{\infty}p^{(m)}$ together with a lift $(p_{l}^{(m)})_{l\in \N} \in \prod A$ of $p^{(m)} \in A_{\infty}$, we apply Lemma \ref{theta-approximation} to obtain a diagram as in (\ref{large-diagram}). Set
\[
h_{l}^{(m)}:= \varrho_{l}^{(m)} \psi_{l}^{(m)}(\be_{Y}), \, l \in \N \, .
\]
This yields a sequence of projections 
\[
p_{l}^{(m+1)}:= p_{l}^{(m)} - h_{l}^{(m)} = \be_{A} - \sum_{k=0}^{m} h_{l}^{(k)} \in A \, ;
\]
let $p^{(m+1)}$ denote the class of $(p_{l}^{(m+1)})_{l\in \N}$ in $A_{\infty}$. By Lemma \ref{theta-approximation} we have 
\[
q((h_{l}^{(m)})_{l\in \N}) \in \theta^{(m)}(\Ch(Y))'\cap p^{(m)}A_{\infty}p^{(m)}
\]
and, since $\theta^{(m)}(\be_{Y})= p^{(m)}$, 
\[
p^{(m+1)} = p^{(m)} - q((h_{l}^{(m)})_{l\in \N}) \in \theta^{(m)}(\Ch(Y))'\cap p^{(m)}A_{\infty}p^{(m)} \, .
\]
We may thus define a $*$-homomorphism 
\[
\theta^{(m+1)}: \Ch(Y) \to p^{(m+1)}A_{\infty}p^{(m+1)} \subset A_{\infty}
\]
by 
\[
\theta^{(m+1)}(\, . \,) := p^{(m+1)}\theta^{(m)}(\, . \,) \, .
\]
Induction yields a sequence of diagrams as in (\ref{large-diagram}), with $*$-homomorphisms 
\[
\theta^{(m)}: \Ch(Y) \to p^{(m)}A_{\infty}p^{(m)} \, ,
\]
zero-dimensional spaces $X^{(m)}$ and $*$-homomorphisms 
\[
\varrho^{(m)}  \kappa^{(m)}: \Ch(X^{(m)}) \to p^{(m)}A_{\infty}p^{(m)} 
\]
for $m \in \N$. Moreover, by \eqref{theta-approximation-estimate} we have 
\[
\liminf_{l \to \infty} \tau(h_{l}^{(m)}) - \mu \cdot \tau(p_{l}^{(m)}) \ge 0 \; \forall \, \tau \in \TT (A), \, m \in \N \, . 
\]
We show by induction over $m$ that
\begin{equation}
\label{trace-estimate-1}
\liminf_{l \to \infty} \tau(\sum_{k=0}^{m} h _{l}^{(k)}) \ge \mu  \cdot \sum_{k=0}^{m} (1-\mu)^{k} \; \forall \, \tau \in \TT (A), \, m \in \N \, .
\end{equation}
For $m=0$ this is true, since by (\ref{theta-approximation-estimate}), 
\[
\liminf_{l \to \infty} \tau(h_{l}^{(0)}) -  \mu \cdot \tau(p_{l}^{(0)}) = \liminf_{l \to \infty} \tau(h_{l}^{(0)}) - \mu \ge 0 \; \forall \, l \in \N, \tau \in \TT (A) \, . 
\]
Suppose now (\ref{trace-estimate-1}) has been verified for some $m \in \N$, then 
\begin{eqnarray*}
\liminf_{l \to \infty} \tau(\sum_{k=0}^{m+1} h_{l}^{(k)}) & = & \liminf_{l \to \infty} \left(\tau(\sum_{k=0}^{m} h_{l}^{(k)}) + \mu \cdot \tau(p_{l}^{(m+1)}) \right. \\
& &\left. + \tau(h_{l}^{(m+1)}) - \mu \cdot \tau(p_{l}^{m+1})  \right)\\
& \stackrel{(\ref{theta-approximation-estimate})}{\ge} & \liminf_{l \to \infty} \left( \tau(\sum_{k=0}^{m} h_{l}^{(k)}) + \mu \cdot \tau(p_{l}^{(m+1)}) \right)\\
& = & \liminf_{l \to \infty} \left( \tau(\sum_{k=0}^{m} h_{l}^{(k)}) + \mu \cdot (\tau(\be_{A} - \sum_{k=0}^{m} h_{l}^{(k)})) \right) \\
& = & \liminf_{l \to \infty} \left( \mu + (1 - \mu) \cdot \tau(\sum_{k=0}^{m} h_{l}^{(k)}) \right) \\
& \stackrel{(\ref{trace-estimate-1})}{\ge} & \mu + \mu \cdot \sum_{k=1}^{m+1} (1-\mu)^{k} \\
& = & \mu \cdot \sum_{k=0}^{m+1} (1- \mu)^{k} \; \forall \, \tau \in \TT (A), \, l\in \N \, .
\end{eqnarray*}
This proves (\ref{trace-estimate-1}) for $m+1$ in place of $m$, hence for all $m \in \N$ by induction. \\
Since $\lim_{m\to \infty} \mu \cdot \sum_{k=0}^{m}(1-\mu)^{k} = 1$, by Dini's theorem there are $l_{0},m_{0} \in \N$ such that 
\begin{equation}
\label{9}
\tau(\sum_{m=0}^{m_{0}} h_{l}^{(m)}) > 1 - \varepsilon \; \forall \, \tau \in 
\TT (A), \, l \ge l_{0} \, .
\end{equation}
Note that by \eqref{large-diagram}, 
\[
q((h_{l}^{(m)})_{l\in \N}) = \varrho^{(m)}  \kappa^{(m)}(\be_{X^{(m)}}) \in A_{\infty} \; \forall \, m \in \N \, .
\]
Set $\widetilde{X} := \coprod_{m=0}^{m_{0}} X^{(m)} $, then 
\[
\Ch(\widetilde{X}) \cong \bigoplus_{m=0}^{m_{0}} \Ch(X^{(m)})
\]
is a commutative $\AF$ algebra, and so is its image under  $\bigoplus_{m=0}^{m_{0}} \varrho^{(m)} \circ \kappa^{(m)}$ in $A_{\infty}$; this image is of the form $\Ch(X) \subset A_{\infty}$ for some zero-dimensional compact space $X$. Let
\[
p_{l}:= \sum_{m=0}^{m_{0}} h_{l}^{(m)} , \, l \in \N \, ,
\]
then $(p_{l})_{l\in \N}$ lifts $p := \be_{X} \in A_{\infty}$. By \ref{theta-approximation} we have $\varrho^{(m)}  \kappa^{(m)}(\be_{X^{(m)}}) \in \theta(\Ch(Y))' \cap A_{\infty} \; \forall \, m=0, \ldots,m_{0}$, whence $p \in \theta(\Ch(Y))' \cap A_{\infty}$, so assertion (i) of the proposition holds. Furthermore, 
\[
p \theta(\Ch(Y)) = \bigoplus_{m=0}^{m_{0}} q((h_{l}^{(m)})_{l \in \N}) \theta(\Ch(Y)) \subset \bigoplus_{m=0}^{m_{0}} \varrho^{(m)}  \kappa^{(m)}(\Ch(X^{(m)})) = \Ch(X) \, ,
\]
whence (ii) holds; we have already verified (iii) in \eqref{9}.  
\end{nproof}
\en

\bn
We are finally prepared to prove the main result of this section.

\begin{nproof} (of Lemma \ref{commutative-tr0-embedding})
We may clearly assume the elements of $\Fh$ to be positive and normalized. Even more, we may assume that $\sum_{a \in \Fh} a = \be_{Y}$. Set $n:= \card \Fh -1$ and let $\Delta^{n} \subset \R^{n+1}$ denote the full $n$-simplex.  Since $\Ch(\Delta^{n})$ is the universal unital commutative $C^{*}$-algebra generated by $n+1$ positive elements adding up to the unit, $C^{*}(\Fh) \subset \Ch(Y)$ is a quotient of $\Ch(\Delta^{n})$. We have $C^{*}(\Fh) \cong \Ch(Y')$ for some compact subspace $Y'$ of $\Delta^{n}$; by \cite{HW}, Theorem III.1 (cf.\ also \cite{KW}, Proposition 3.3 and \cite{Wi1}, Proposition 2.19), we see that 
\[
\dim Y' \le \dim \Delta^{n} = n \, .
\]
Restricting $\theta$ to $\Ch(Y')$ we see that it suffices to prove the assertion of the lemma for a finite-dimensional space; thus we assume that $\dim Y = n $ for some $n \in \N$ right away. \\
Let $\Ch(X) \subset A_{\infty}$ be a commutative $\AF$ algebra and $(p_{l})_{l\in \N} \subset A$ a sequence of projections satisfying the assertions of Proposition \ref{cantor-interpolation}. By \ref{cantor-interpolation}(i) we then have $\be_{X} \in \theta(\Ch(Y))' \cap A_{\infty}$; by (ii), $\be_{X} \theta(\Ch(Y)) \subset \Ch(X)$. But then there is a unital finite-dimensional $C^{*}$-subalgebra $\bar{B}$ of $\Ch(X)$ such that $\dist(\be_{X} \theta(a) \be_{X},\bar{B}) < \varepsilon/2 \; \forall \, a \in \Fh$; note that $\bar{B}$ is commutative and that we may assume $\be_{\bar{B}} = \be_{X}$. Since $(p_{l})_{l \in \N}$ lifts $\be_{X}$, $\be_{X}A_{\infty}\be_{X}$ is a quotient of $\prod_{l \in \N} p_{l}Ap_{l}$. Finite-dimensional $C^{*}$-algebras are semiprojective (\cite{Lo}, Chapter 14), so there is a $*$-homomorphism 
\[
\sigma: \bar{B} \to \prod_{l \in \N} p_{l}Ap_{l}
\]
lifting $\id_{\bar{B}}$. It is clear that all but finitely many of the $\sigma_{l}$ (the components of $\sigma$) are unital. Now if $l_{0} \in \N$ is large enough, then $\dist(\theta(a),\sigma_{l_{0}}(\bar{B}))< \varepsilon \; \forall \, a \in \Fh$, $\sigma_{l_{0}}(\be_{\bar{B}}) = p_{l_{0}}$ and $\|[\theta(a),p_{l_{0}}]\|< \varepsilon \; \forall \, a \in \Fh$. Since $\tau(p_{l_{0}})>1-\varepsilon \; \forall \, \tau \in \TT (A)$ (for large enough $l_{0}$) by Proposition \ref{cantor-interpolation}(iii), $B:= \sigma_{l_{0}}(\bar{B})$ satisfies the assertion of the lemma.    
\end{nproof}
\en

\section{The uniqueness theorem}

    Towards the uniqueness theorem, we need the following stable
uniqueness result which can be found in \cite{GongLin2} Theorem 3.1 (see 
also \cite{GongLin2} Remark 1.1, 
\cite{GongLin} and \cite{Linhomomorphism}):

\bn
\begin{thm}  Let $X$ be a compact metric space.  For any $\epsilon >0$
and any finite subset $\F \subseteq C(X)$, there exist $\delta > 0$, 
$\eta > 0$, an integer $N > 0$, a finite subset $\G \subseteq C(X)$
and a finite subset $\P \subseteq \mathbb{P}(C(X))$ satisfying the 
following:

   For any unital simple separable nuclear
$C^*$-algebra $\A$ with real rank zero, stable rank one
and weakly unperforated $\Kk_0$ group,
for any $\eta$-dense subset $\{ x_1, x_2, .., x_k \}$ in $X$,   
and any $\G - \delta$-multiplicative
c.p.c. linear maps $\phi, \psi : C(X) \rightarrow
\A$, if 
\[
[ \phi ] | \P = [ \psi ] | \P,
\]
then there exists a unitary $u \in M_{Nk + 1} ( \A )$ such that
\[
u( \phi(f) \oplus f(x_1) 1_N \oplus f(x_2) 1_N \oplus ... f(x_k) 1_N) u^*
\approx_{\epsilon}  
\psi(f) \oplus f(x_1) 1_N \oplus f(x_2) 1_N \oplus ... f(x_k) 1_N
\]
for all $f \in \F$.  
\label{thm:Linstableuniqueness}  
\end{thm}
\en

\bn
\label{d-tr0-embedding}
\begin{df}  Let $X$ be a compact metric space, and
let $\A$ be a unital simple separable $C^*$-algebra.
A $*$-monomorphism $\phi : C(X) \rightarrow \A$ is said to be a  
\emph{$\TAF$-embedding} (i.e., a \emph{tracially $\AF$-embedding}) if 
for every $\epsilon > 0$, for every finite subset $\F \subseteq \A$ and
for every nonzero positive element $a \in \A$, 
there is a projection $p \in \A$ and there is a finite dimensional
$C^*$-subalgebra $F \subseteq \A$ such that 
\begin{enumerate}
\item[(i)] $p$ is Murray-von Neumann equivalent to a subprojection of 
$\her(a)$ (the hereditary subalgebra of $\A$ generated by $a$),
\item[(ii)] $\be_{F} = \be - p$,
\item[(iii)] $\| pa - ap \| < \epsilon$ for all $a \in \F$, and 
\item[(iv)] $pap$ is within $\epsilon$ of an element of $F$ for all $a \in \F$.
\end{enumerate}
\end{df}

    From Lemma \ref{commutative-tr0-embedding}, it follows that 
if $A$ is a unital separable simple real rank zero $\Zh$-stable 
$C^*$-algebra, and if $X$ is a compact second countable metric space, then
any unital $*$-monomorphism $\phi : C(X) \rightarrow A$ is a
$\TAF$-embedding.
\en

\bn
\begin{prop}
Let $X$ be a compact metric space, and let $\A$ be a unital
separable simple nuclear
$C^*$-algebra with real rank zero, stable rank one and
weakly unperforated $\Kk_0$ group.
Let $\phi : C(X) \rightarrow \A$ be a unital 
$\TAF$-embedding.\\  
Then for every $\epsilon > 0$, for every finite subset $\F \subseteq \A$,
for every integer $N \geq 1$ and for every nonzero positive element 
$a \in \A$, there exists a real number
$\delta$ with $ 0 < \delta < \epsilon$, there exists a finite subset
$\{ x_1, x_2, ..., x_n \} \subseteq X$ which is $\delta$-dense in $X$, and
there exists a projection $p \in \A$ and a $\F - \epsilon$-multiplicative
c.p.c. $L_1 : C(X) \rightarrow p \A p$
such that 
\begin{enumerate} 
\item[(i)] $p$ is Murray-von Neumann equivalent to a subprojection of 
$\her(a)$ (there hereditary subalgebra of $\A$ generated by $a$),
\item[(ii)] $\| p \phi(f) - \phi(f)p \| < \epsilon$ for all $f \in \F$, and
\item[(iii)] there exists pairwise orthogonal projections 
$p_0, p_1, p_2, ..., p_n, t$ with 
\[
p_0 = p, \makebox{ }\sum_{i=0}^n p_i + t = \be_{\A},
\makebox{ and }
Np \preceq p_i
\]
 for $i \neq 0$ (here $\preceq$ is the relation of 
being Murray-von Neumann equivalent to a subprojection) and there exists
a finite dimensional $*$-homomorphism 
\[
h_1 : C(X) \rightarrow t \A t
\]
such that\\ 
$\phi(f)$ is within $\epsilon$ of 
$L_1 (f) + \sum_{i=1}^n f(x_i)p_i + h_1(f)$, 
for all $f \in \F$.
\end{enumerate}
\label{prop:TAFmorphisms}
\end{prop}
\begin{nproof} 
For simplicity, let us assume that the elements of $\F$ all have norm 
less than or equal to one.
Let $d$ be the metric on the space $X$.
Let $\delta$ be a real number such that  $0 < \delta < \epsilon/10$ and 
such that for every $x,y \in X$, if $d(x,y) < \delta$ then 
$| f(x) - f(y) | < \epsilon/10$.  
Now let $\{ x_1, x_2, ..., x_n \}$ be a $\delta$-dense subset of $X$.
We may assume that $d(x_i, x_j) >0$ for $i \neq j$.  
For each $i$, let $0 < \delta_i < \delta$ be such that 
$\overline{B(x_i, \delta_i)} \cap \overline{ B(x_j, \delta_j)} = 
\emptyset$ for $i \neq j$. 
(Here, $\overline{B(x_i, \delta_i)}$ is the closure of the open ball
$B(x_i, \delta_i)$
with radius $\delta_i$ about $x_i$.)
For each $i$, let $f_i$ be a positive real-valued function with
$0 \leq f_i \leq 1$, $f_i(x_i) = 1$ and $supp(f_i) \subseteq B(x_i, \delta_i)$,
where $supp(f_i)$ is the (compact subset of $X$ which is the) support
of $f_i$.
Since $\phi$ is a $*$-monomorphism and since $\A$ is simple, 
for each $i$, $inf\{ \tau( \phi( f_i )) : \tau \in \TT (\A) \} > 0$ (strictly
greater than zero).
Hence, let $s$ be the strictly positive real number given by
$s =_{df} inf \{ \tau( \phi(f_i)) : \tau \in \TT (\A), 1 \leq i \leq n \}$.
Let $\F' =_{df} \F \cup \{ f_i : 1 \leq i \leq n \}$.

    Now apply Lemma \ref{lem:LinbookMatrix} on $X$, $\F'$ and
$\epsilon_1 =_{df} \mmin\{ \epsilon/10, s/(10(N + 10)) \}$ to get
$\G$ and $\rho$ ($\rho$ is the $\delta$ in 
Lemma $\ref{lem:LinbookMatrix}$).  We may assume that $\rho < \epsilon_1$. 
Making $\rho$ smaller if necessary, we may assume that the elements of 
$\G$ all have norm less than or equal to one.
Since the $*$-monomorphism $\phi : C(X) \rightarrow \A$ is a  
$\TAF$-embedding and since finite dimensional $C^*$-algebras are 
injective von Neumann algebras, let $p \in \A$ be a projection and let
$F \subseteq \A$ be a finite dimensional $C^*$-subalgebra with 
$\be_{F} = \be - p$ such that 
\begin{enumerate}
\item[i.] $p$ is Murray-von Neumann equivalent to a subprojection of 
$\her(a)$, the hereditary subalgebra of $\A$ generated by $a$,
\item[ii.] $\tau(p) < \epsilon_1$ for all $\tau \in \TT (\A)$,
\item[iii.] $\| p \phi(f) - \phi(f) p \| < \epsilon_1$ 
for all $f \in \G$,
and 
\item[iv.]  there exists a unital c.p.c. 
$\G - \rho$-multiplicative map $L : C(X) \rightarrow F$ such that 
$\phi(f)$ is  
within $\epsilon_1$ of $p \phi(f) p + L(f)$ for all $f \in \G$.
\end{enumerate}

   Applying Lemma \ref{lem:LinbookMatrix}  to $L$, we get a $*$-homomorphism
 $h : C(X) \rightarrow F$ with a projection $q = h(\be_{C(X)})$ such that
$\tau'(\be - p - q) < \epsilon_1$ for all $\tau' \in \TT (F)$ and 
such that  
$L(f)$ is within $\epsilon_1$ of 
$(\be - p - q) L(f) (\be - p - q) + h(f)$ for all $f \in \F'$.   
Note that we must have that $\tau(\be - p - q) < \epsilon_1$ for all
$\tau \in \TT (\A)$.   

    Hence, $\phi(f)$ is within $2 \epsilon_1$ 
of $L'(f) + h(f)$ for all $f \in \F'$, where 
$$L'(f) = _{df} p \phi(f) p +  
(\be - p - q) L(f) ( \be - p - q)$$ (so $L'(\be_{C(X)}) = \be - q$).   
Note that (after some computation)
$L'$ is $\F - 3\epsilon_1$-multiplicative.
Also, $\tau(q) > (1 - \epsilon_1)^2$ (which in turn is greater than
$1 - 2 \epsilon_1$) 
for all $\tau \in \TT (\A)$.     

Now for all $i$, 
$\phi(f_i)$ is within $2 \epsilon_1$ of $L'(f_i) + h(f_i)$.  Hence,
for all $i$ and for all $\tau \in \TT (\A)$, 
$\tau( \phi(f_i))$ is within $2 \epsilon_1$ of $\tau(L'(f_i)) + 
\tau(h(f_i))$.
Also, note that $\tau(L'(f_i)) < 1 - (1 - \epsilon_1)^2 = 2 \epsilon_1 - 
{\epsilon_1}^2$.      
Hence, for all $i$ and for all $\tau \in \TT (\A)$, 
\begin{eqnarray*}
\tau(h(f_i)) & = & ( \tau(h(f_i)) + \tau(L'(f_i)) ) - \tau(L'(f_i)) \\
            &  \geq & (\tau(h(f_i)) + \tau(L'(f_i)) - 
2\epsilon_1 + {\epsilon_1}^2 \\
             & \geq & (\tau(\phi(f_i)) - 2 \epsilon_1) - 2 \epsilon_1 +
{\epsilon_1}^2 \\
             & = & \tau(\phi(f_i)) - 4\epsilon_1 + {\epsilon_1}^2 \\      
             & \geq & s - 4\epsilon_1 + {\epsilon_1}^2  
\end{eqnarray*} 
Recall that $\epsilon_1 = \mmin\{ \epsilon/10, s/(10(N + 10)) \}$.
Hence, for all $\tau \in \TT (\A)$ and for all $i$, 
\begin{eqnarray*}
\tau(h(f_i)) & \geq & s - 4s/(5(N + 10))  \\
             & \geq & (5Ns + 50s - 4s)/ (5(N + 10)) \\
             & = &  (5N + 46)s/ (5N + 50) 
\end{eqnarray*}
             
     Now suppose that $r_1, r_2, ..., r_l$ are pairwise orthogonal projections
in
$h(\be_{C(X)}) F h ( \\ \be_{C(X)})$ and 
$y_1, y_2, ..., y_l$ points in $X$ such 
that for all $f \in C(X)$, 
$h(f) = \sum_{j=1}^l f(y_j) r_j$.
For each $i$, let $S_i =_{df} \{ j : d(y_j, x_i) < \delta_i \}$.
Then for each $i$ and for each $\tau \in \TT (\A)$,
$\tau( \sum_{j \in S_i } r_j ) \geq (5N + 46)s/ (5N + 50)$.    
Replace $h$ by a finite-dimensional $*$-homomorphism $h' : C(X) \rightarrow
F$ where for all $f \in C(X)$, 
$$h'(f) =_{df} \sum_{i=1}^n f(x_i) \sum_{j \in S_i} r_j + \sum_{j \not\in 
S_i \forall i} f(y_j) r_j$$ 
Hence, for each $f \in \F$, $h(f)$ is within $\epsilon/10$ of 
$h'(f)$.  
Recall that $L'(\be) = \be - q$ and $\tau(q) > 1- 2\epsilon_1 \geq
1 - s/(5N + 50)$ for all $\tau \in \TT (\A)$. 
Hence, $\tau(L'(\be)) < s/(5N + 50)$ for all $\tau \in \TT (\A)$.   
Hence, $N \tau(L'(\be)) < \tau( \sum_{j \in S_i } r_j )$ for all 
$\tau \in \TT (\A)$.  Hence, 
for each $f \in \F$,
$\phi(f)$ is within $3 \epsilon/10$ of $L'(f) + h'(f)$ and the approximation 
$L' + h'$ is what is required. 
\end{nproof}

\en

\bn
\begin{nproof}(of Theorem \ref{thm:MainUniqueness})
    The ``only if" direction follows from \cite{Rorpi} 5.4.
Hence, we need only prove the ``if" direction. 

    Firstly, by \cite{Linembedding}
Theorem 4.5, there is a unital simple $\AH$-algebra
$\A_0$ with bounded dimension growth and real rank zero and there is 
a unital $*$-homomorphism $\Phi : \A_0 \rightarrow \A$ such that
$\Phi$ induces an order isomorphism (which, of course, respects the Bockstein
operations) between the full $\mathrm{K}$-groups $\underline{\Kk}(\A_0)$ and
$\underline{\Kk}(\A)$; and since $\A_0$ and $\A$ are both real rank zero,
this induces an isomorphism between the tracial simplexes $\TT (\A)$ and 
$\TT (\A_0)$.  Henceforth, to simplify notation, we may identify 
$\underline{\Kk}(\A_0)$, $\TT (\A_0)$ with $\underline{\Kk}(\A)$, $\TT (\A)$ 
respectively and take the induced maps  $\underline{\Kk}(\Phi)$, $\TT (\Phi)$
to be identity maps.  
Hence, by Theorem \ref{thm:Existence}, 
let $\theta : C(X) \rightarrow
\A_0 \subseteq \A$ be a unital $*$-monomorphism such that 
$\KL (\theta) = \KL (\phi)$ ($= \KL (\psi)$) and 
$\tau \circ \theta = \tau \circ \phi$ ($= \tau \circ \psi$) for all 
$\tau \in \TT (\A)$.   
To prove that $\phi$ and $\psi$ are approximately unitarily equivalent,
it suffices to show that $\phi$ and $\theta$ are approximately 
unitarily equivalent. 

    Let $\epsilon > 0$ and a finite subset $\F \subseteq C(X)$ be given.
Let us assume that the elements of $\F$ all have norm less than or equal
to one.
Let $\nu > 0$ be such that $| f(x) - f(y) | < \epsilon/8$ if 
$\dist (x,y) < \nu$ for all $x,y \in X$ and all $f \in \F$.
We may assume that $\nu < 1$
Now
\begin{enumerate}
\item[i.] take $s = l$; 
\item[ii.] take a finite $\nu/2$-dense set $\{ x_1, x_2, ..., x_m \}$ in $X$
for which $O_i \cap O_j = \emptyset$ whenever $i \neq j$, where
$O_i =_{df} \{ x \in X : \dist (x, x_i) < \nu/2 \}$;  and also, put
$\tilde{O_i} =_{df} \{ x \in X : \dist (x, x_i) < \nu/4 \}$.  
\item[iii.] and take a real number $\sigma$ where 
$0 < \sigma < 1/2$ and where
$\mu_{\tau \circ \phi}( \tilde{O_i}) > 10\sigma$ 
for all $\tau \in \TT (\A)$ and for
all $i$.  
\end{enumerate}
   Taking the above data and putting them into Theorem 
\ref{thm:Linapproximateuniqueness},    
we get quantities
$\gamma$, $\F_1$, $\epsilon_1$ and $\P_1$
($\F_1$, $\epsilon_1$, $\P_1$ are the 
$\G$, $\delta$, $\P$ respectively in the conclusion of 
Theorem \ref{thm:Linapproximateuniqueness}).   
For each $i$, let $f_i$ be the real-valued function on $X$ with
(a) $0 \leq f_i \leq 1$, (b) $f_i = 1$ on $\tilde{O_i}$ and
(c) the support of $f_i$ is contained in $O_i$.  
Expanding $\F_1$ if necessary, we may assume that $\F_1$ contains
$\{ f_i : 1 \leq i \leq m \} \cup \F$.  We may also assume that all the 
elements of $\F_1$ have norm less than or equal to one.

     Now put $\epsilon$ and  $\F_1$ into Theorem 
\ref{thm:Linapproximateuniqueness} 
to get quantities
$\epsilon_2$, $\eta$, $N > 0$, $\F_2$ and $\P_2$ ($\epsilon_2$, $\F_2$ 
$\P_2$ are the
$\delta$, $\G$ $\P$ 
in the conclusion of Theorem \ref{thm:Linapproximateuniqueness}).
Contracting $\epsilon_2$ if necessary, we may assume that 
$\epsilon_2 < \epsilon_1$.   Expanding $\F_2$ if necessary,
we may assume that $\F_1 \subseteq \F_2$ and that the elements of
$\F_2$ all have norm less than or equal to one.
Finally, expanding $\P_2$ if necessary, we may assume that $\P_1 \subseteq
\P_2$.

    Let $\epsilon_3 > 0$ and let $\F_3$ be a finite subset of $C(X)$ 
such that
\begin{enumerate}
\item $\epsilon_3 < \epsilon_2$,
\item $\epsilon_3 < \mmin \{ \gamma/100, \sigma/100, \eta/100 \}$, 
\item $\F_2 \subseteq \F_3$, and
\item if $\rho_1, \rho_2 : C(X) \rightarrow \A$ are two 
$\F_3 - \epsilon_3$-multiplicative c.p.c. maps
then $\rho_1$ and $\rho_2$ are both well-defined on $\P_2$.
Moreover, if $\rho_1(f)$ is within $\epsilon_3$ of $\rho_2(f)$ for all
$f \in \F_3$ then
\[
[ \rho_1] | \P_2 = [ \rho_2] | \P_2.
\]
\end{enumerate}
   
   Now by Lemma \ref{commutative-tr0-embedding}, $\phi$ is a $\TAF$-embedding. 
Therefore,  put $\phi$,
$\epsilon_3$, $\F_3$, $N$,  and a projection whose value at
every unital trace is strictly less than $\epsilon_3$ into Proposition 
\ref{prop:TAFmorphisms}.
Then as a consequence, we get a real number $\delta$ with 
$0 < \delta < \epsilon_3$, a $\delta$-dense subset $\{ x_1, x_2, ..., x_n \}$
of $X$, and  we get a projection $p \in \A$ and a 
$\F_3 - \epsilon_3$-multiplicative c.p.c. 
map $L_1 : C(X) \rightarrow p\A p$ such that
\begin{enumerate}
\item[(a)] $\tau(p) < \epsilon_3$ for all $\tau \in \TT (\A)$,
\item[(b)] there exists pairwise orthogonal projections $p_0, p_1, p_2, ...,
p_n, t$ with $p_0 = p$, $\sum_{i=0}^n p_i + t = \be_{\A}$, and 
$Np \preceq p_i$ for $i \neq 0$ (here $\preceq$ is the relation of being
Murray-von Neumann equivalent to a subprojection) and there exists a unital 
finite dimensional $*$-homomorphism $h_1 : C(X) \rightarrow t \A t$
such that 
\[
\phi(f) \approx_{\epsilon_3} L_1(f) + \sum_{i=1}^n f(x_i) p_i + h_1(f)
\]
for all $f \in \F_3$.  
\end{enumerate}
 
   Let $\phi_1 : C(X) \rightarrow \A$ be the c.p.c. 
$\F_3 - \epsilon_3$-multiplicative map given by 
$\phi_1 : f \mapsto L_1(f) + \sum_{i=1}^n f(x_i) p_i + h_1(f)$.
By the definition of $\epsilon_3$, $\phi_1$ (and also $L_1$) is 
well-defined on $\P_3$.  Moreover,
\[
[\phi] | \P_3 = [ \phi_1 ] | \P_3.
\]
Also, since $\tau \circ \phi(f_i) > 10 \sigma$ for all $\tau \in \TT (\A)$
and for all $i$, we must have that $\tau \circ \phi_1(f_i) > 9 \sigma$
for all $\tau \in \TT (\A)$ and for all $i$.

     Now let $q_0, q_1, q_2, ..., q_n, t'$ be pairwise 
orthogonal projections
in $\A_0$ such that (a) $\be_{\A_0} = \sum_{i=1}^n q_i + t'$,  
(b) $q_i$ is Murray-von Neumann equivalent (in $\A$) to $p_i$, and
(c) $t'$ is Murray-von Neumann equivalent (in $\A$)  to $t$.

     Since $X$ is a path connected,  
$\KL (\phi) - \KL (h_1) \in \KL ( C(X), q_0 \A_0 q_0 )_{+,1}$.  
Hence, by Theorem \ref{thm:Existence},  let 
$L_2 : C(X) \rightarrow q_0 \A q_0$ be a unital $*$-homomorphism
such that $\KL (L_2) = \KL (\phi) - \KL (h_1)$.
Let $h_2 : C(X) \rightarrow t' \A t'$ be a unital finite dimensional 
$*$-homomorphism which is unitarily equivalent (with unitary in $\A$)
to $h_1$.
Now consider the $*$-homomorphism 
$\phi_2 : C(X) \rightarrow \A$ given by 
$\phi_2 : f \mapsto L_2(f) + \sum_{i=1}^n f(x_i) q_i + h_2(f)$.
Clearly, 
\[
[\phi_1 ] | \P_3 = [ \phi_2 ] | \P_3.
\]
Hence, by Theorem \ref{thm:Linstableuniqueness}, 
let $u \in \A$ be a unitary
such that 
\[
u \phi_1(f) u^* \approx_{\epsilon_2} \phi_2(f)
\]
for all $f \in \F_3$.
Hence, since $\tau \circ \phi_1 (f) > 9 \sigma$ for all $f \in \F_3$ and 
for all $\tau \in \TT (\A)$,
$\phi_2(f) > 8 \sigma$ for all $f \in \F_3$ and for all $\tau \in \TT (\A)$.
Hence, $\mu_{\tau \circ \phi_2} (O_i) > \sigma \nu$ for all $f \in F_3$, 
all $\tau \in \TT (\A)$ and all $i$. 

      Also,
\begin{eqnarray*}
\lefteqn{\| \tau \circ \theta (f) - \tau \circ \phi_2(f) \|}  \\
& = & \| \tau \circ \phi(f) - \tau \circ \phi_2(f) \|   \\ 
& \leq & \| \tau \circ \phi(f) - \tau \circ \phi_1(f) \| +  
\| \tau \circ \phi_1(f) - \tau \circ \phi_2(f) \|  \\
& < & \epsilon_3 + \| \tau \circ \phi_1(f) - \tau \circ \phi_2 (f) \|  \\
& < & \epsilon_3  + \epsilon_2 \\  
& \leq & \gamma,  
\end{eqnarray*}
for all $f \in \F_3$.

    Finally,   
\[
[\theta] | \P_3 = [ \phi ]| \P_3 = [ \phi_1 ] | \P_3 = [ \phi_2 ] | \P_3.
\]

     Hence, it follows, by Theorem \ref{thm:Linapproximateuniqueness},
that there is a unitary $v \in \A_0$ such that 
\[ 
 v \phi_2(f) v^* \approx_{\epsilon} \theta(f),
\]
for all $f \in \F$.
Hence, 
\[
v u \phi(f) u^* v^* \approx_{\epsilon} \theta(f),
\]
for all $f \in \F$.
Since $\epsilon$ anf $\F$ were arbitrary,
$\phi$ and $\theta$ are approximately unitarily equivalent.
Similar for $\psi$ and $\theta$.
Hence $\phi$ and $\psi$ are approximately unitarily equivalent as required.
\end{nproof}

\en

\bn
   Finally, note that since a finite CW complex is the finite union of 
pairwise disjoint, path connected, second countable, compact metric spaces,
we can replace $X$ in Theorem \ref{thm:MainUniqueness} by an arbitrary 
finite CW complex. 




\en

\end{document}